\documentclass[a4paper,draft,12pt,reqno]{article}
\usepackage{amsmath,amssymb,amsthm,xspace,a4wide,amscd,latexsym,
amssymb}
\usepackage[all]{xy}
\usepackage{enumerate}

\DeclareMathOperator{\Hom}{Hom}

\DeclareMathOperator{\Rad}{Rad}
\renewcommand{\hom}{\operatorname{hom}}
\DeclareMathOperator{\Soc}{Soc}

\DeclareMathOperator{\soc}{soc}

\newcommand{\tto}{\twoheadrightarrow}

\newcommand{\Oo}{\ensuremath{\mathcal{O}}}

\newcommand{\ppp}{\ensuremath{\mathfrak{p}}}

\newcommand{\Z}{\ensuremath{\mathbb{Z}}}

\renewcommand{\phi}{\varphi}

\newcommand{\End}{\operatorname{End}}

\newcommand{\Ext}{\operatorname{Ext}}
\newcommand{\ext}{\operatorname{ext}}
\newcommand{\rank}{\operatorname{rank}}

\begin{document} 
\theoremstyle{plain}
\newtheorem{example}{Example}
\newtheorem{definition}{Definition}
\newtheorem{remark}{Remark}
\newtheorem{theorem}{Theorem}
\newtheorem{lemma}{Lemma}
\newtheorem{proposition}{Proposition}
\newtheorem{corollary}{Corollary}
\newtheorem{conjecture}{Conjecture}

\title{A pairing in homology and the category of linear complexes 
of tilting modules for a quasi-hereditary algebra}
\author{Volodymyr Mazorchuk and Serge Ovsienko}
\date{}
\maketitle

\begin{abstract}
We show that there exists a natural non-degenerate pairing of 
the homomorphism space between two neighbor standard modules over a 
quasi-hereditary algebra with the first extension space between the
corresponding costandard modules and vise versa. Investigation
of this phenomenon leads to a family of pairings involving 
standard, costandard and tilting modules. In the graded case,
under some "Koszul-like" assumptions (which we prove are 
satisfied for example for the blocks of the category $\mathcal{O}$),
we obtain a non-degenerate pairing between certain graded 
homomorphism and graded extension spaces. This motivates the study
of the category of linear tilting complexes for graded 
quasi-hereditary algebras. We show that, under assumptions,
similar to those mentioned above, this category realizes
the module category for the Koszul dual of the Ringel dual 
of the original algebra. As a corollary we obtain that under
these assumptions the Ringel and Koszul dualities commute.
\end{abstract}

\section{Introduction and description of the results}\label{s1} 

Let $\Bbbk$ be an algebraically closed field. If the opposite is 
not emphasized, in this paper by a module we mean a {\em left} module 
and we denote by  $\Rad(M)$ the radical of a module, $M$. For a 
${\Bbbk}$-vector space,  $V$, we denote the dual space by $V^*$.

Let $A$ be a basic $\Bbbk$-algebra, which is  quasi-hereditary 
with respect to the natural order on the indexing set 
$\{1,2,\dots,n\}$ of pairwise-orthogonal  primitive idempotents 
$e_i$ (see \cite{CPS1,DR,DR2} for details). Let $P(i)$, $\Delta(i)$, 
$\nabla(i)$, $L(i)$, and $T(i)$ denote the projective, standard, 
costandard, simple and tilting $A$-modules, associated to $e_i$, 
$i=1,\dots, n$, respectively. Set $P=\oplus_{i=1}^n P(i)$,
$\Delta=\oplus_{i=1}^n \Delta(i)$,$\nabla=\oplus_{i=1}^n \nabla(i)$,
$L=\oplus_{i=1}^n L(i)$, $T=\oplus_{i=1}^n T(i)$.
We remark that, even if the standard $A$-modules are fixed, the
linear order on the indexing set of primitive idempotents, with respect
to which the algebra $A$ is quasi-hereditary, is not unique in general.
We denote by $R(A)$ and $E(A)$ the Ringel and Koszul duals of $A$ 
respectively. A graded algebra, $B=\oplus_{i\in\Z}B_i$, will be called 
{\em positively graded} provided that $B_i=0$ for all $i<0$ and 
$\Rad(B)=\oplus_{i>0}B_i$.

This paper has started from an attempt to give a conceptual explanation 
for the equality 
\begin{equation}\label{introeq1}
\dim\Hom_{A}(\Delta(i-1),\Delta(i))=
\dim\Ext^1_A(\nabla(i),\nabla(i-1)),
\end{equation}
which is proved at the beginning of Section~\ref{s2}. Our first main 
result, proved also in Section~\ref{s2}, is the following statement:

\begin{theorem}\label{introt1}
\begin{enumerate}[(1)]
\item\label{introt1.1} Let $i,j\in\{1,\dots,n\}$ and $j<i$. Then there 
exists a bilinear pairing, 
\begin{displaymath}
\langle\cdot\, , \cdot\rangle:
\Hom_A\left(\Delta(j),\Delta(i)\right)\times
\Ext_A^1\left(\nabla(i),\nabla(j)\right)\to\Bbbk.
\end{displaymath}
\item\label{introt1.2} If $j=i-1$, then $\langle\cdot\, ,\cdot\rangle$
is non-degenerate. 
\end{enumerate}
\end{theorem}

Theorem~\ref{introt1} explains the origins of \eqref{introeq1} and motivates
the study of $\langle\cdot\, ,\cdot\rangle$. It happens that in the
general case, that is for $j<i-1$, the analogue of 
Theorem~\ref{introt1}\eqref{introt1.2} is no longer true. We give an example at 
the beginning of Section~\ref{s3}. In the same section we present some 
special results and a modification of $\langle\cdot\, ,\cdot\rangle$
in the general case. 

An attempt to lift the above results to higher $\Ext$'s naturally led us
to the definition of a different pairing, which uses a minimal tilting resolution
of the costandard module. In Section~\ref{s4} we construct and investigate
a pairing between $\Ext^l_A(\nabla,\nabla)$ and $\Hom_A(\Delta,T_l)$, where
$T_l$ is the $l$-th component of a minimal tilting resolution of $\nabla$. 
In the case $l=1$ this new pairing induces the one we have constructed
in Section~\ref{s2}. 

The new pairing is rarely non-degenerate. In an attempt to find some 
conditions, which would ensure this property, we naturally came to the 
graded case. In Section~\ref{s5} we show that in the graded case our new 
pairing induces  a non-degenerate pairing  between the graded homomorphism 
and the graded first extension spaces under the condition that the
costandard modules admit linear tilting resolutions.  Here the linearity of 
the resolution means the following: we show that for a positively graded
quasi-hereditary algebra all tilting modules are gradable and thus we can
fix their graded lifts putting their "middles" in degree $0$; the linearity
of the resolution now means that the $i$-th term of the resolution consists
only of tilting modules, whose "middles" are exactly in degree $i$.
This observation brings the linear complexes of tilting modules into
the picture and serves as a bridge to the second part of the paper, in which 
we study the category of all such linear complexes.

The above mentioned condition of the existence of a linear tilting
resolution for costandard $A$-modules immediately resembles the conditions,
which appeared in \cite{ADL} during the study of the following question: when
the Koszul dual of a quasi-hereditary algebra is quasi-hereditary with
respect to the opposite order? In \cite[Theorem~3]{ADL} it was shown that 
this is the case if and only if both the standard and costandard $A$-modules
admit a linear projective and injective (co)resolution respectively
(algebras, satisfying these conditions, were called {\em standard Koszul}
in \cite{ADL}). This resemblance motivated us to take a closer look at the 
category of linear complexes of tilting $A$-modules. The most striking property 
of this category is the fact that it combines two objects of completely different
natures: tilting modules for a quasi-hereditary algebra, which give rise
to the so-called {\em Ringel duality}; and linear resolutions, which are
the source of a completely different duality, namely the {\em Koszul duality}.
Under some natural assumptions, which roughly mean that all objects we consider 
are well-defined and well-coordinated with each other, in 
Section~\ref{s55} we prove our second main result:

\begin{theorem}\label{introt2}
Assume that $A$ is a positively graded quasi-hereditary algebra, such that
\begin{enumerate}[(i)]
\item standard $A$-modules admit a linear tilting coresolution,
\item costandard $A$-modules admit a linear tilting resolution.
\end{enumerate}
The above conditions imply that the quadratic dual $R(A)^!$ of $R(A)$ is 
quasi-hereditary (with respect to the same order as for $A$), and we
further assume that
\begin{enumerate}[(i)]
\setcounter{enumi}{3}
\item the grading on $R(R(A)^!)$, induced from the category of graded 
$R(A)^!$-modules, is positive.
\end{enumerate}
Then the algebras $A$, $R(A)$, $E(A)$, $R(E(A))$ and $E(R(A))$ are standard 
Koszul quasi-hereditary algebras, moreover, $E(R(A))\cong R(E(A))$ as 
quasi-hereditary algebras. In other words, Koszul and Ringel dualities commute 
on $A$.
\end{theorem}

As a preparatory result for this theorem we show that, under the same
assumptions, the category of bounded linear complexes of tilting $A$-modules
is equivalent to the category of graded modules over $E(R(A))^{opp}$.
Moreover, this realization preserves (in some sense) standard and costandard
modules but switches simple and tilting modules.

We finish the paper with proving that all conditions of Theorem~\ref{introt2} 
are satisfied for the associative algebras, associated with the blocks of the 
BGG category $\mathcal{O}$. This is done in Section~\ref{s6}.
In the same section we also derive some consequences for these algebras, 
in particular, about the structure of tilting modules. The paper is finished
with an Appendix, written by Catharina Stroppel, where it is shown that 
all conditions of Theorem~\ref{introt2} are satisfied for the associative 
algebras, associated with the blocks of the parabolic category $\mathcal{O}$
in the sense of \cite{RC}. As the main tool in the proof of the last result, 
it is shown that Arkhipov's twisting functor on $\mathcal{O}$ (see \cite{AS,KM}) 
is gradable.

For an abelian category, $\mathcal{A}$, we denote by
$D^b(\mathcal{A})$ the corresponding bounded derived category and by
$K(\mathcal{A})$ the corresponding homotopic category.
In particular, for an associative algebra, $B$, we denote by $D^b(B)$ 
the bounded derived category of $B\mathrm{-mod}$ and by $K(B)$ the homotopic
category of $B\mathrm{-mod}$. For $M\in B\mathrm{-mod}$ 
we denote by $M^{\bullet}$ the complex defined via
$M^{0}=M$ and $M^{i}=0$, $i\neq 0$.

We will say that a module, $M$, is {\em Ext-injective} (resp. 
{\em Ext-projective}) {\em with respect to a module}, $N$, provided that 
$\Ext^{k}(X,M)=0$, $k>0$, (resp. $\Ext^{k}(M,X)=0$, $k>0$) for any 
subquotient $X$ of $N$. 

When we say that a graded algebra is Koszul, we mean that it is Koszul
with respect to this grading.

\section{A bilinear pairing between $\Hom_A$ and $\Ext^1_A$}\label{s2} 

The following observation is the starting point of this paper.
Fix $1<i\leq n$. According to the classical BGG-reciprocity for 
quasi-hereditary algebras (see for example \cite[Lemma~2.5]{DR2}), 
we have that  $[I(i-1):\nabla(i)]=[\Delta(i):L(i-1)]$,
where the first number is the multiplicity of $\nabla(i)$ in a 
costandard filtration of $I(i-1)$, and the second number is the
usual composition multiplicity. The quasi-heredity of $A$,
in particular, implies that $\Delta(i-1)$ is Ext-projective with
respect to $\Rad(\Delta(i))$ and 
hence $[\Delta(i):L(i-1)]=\dim\Hom_{A}(\Delta(i-1),\Delta(i))$.

The number $[I(i-1):\nabla(i)]$ can also be reinterpreted. Again,
the quasi-heredity of $A$ implies that any non-zero element from
$\Ext^1_A(\nabla(i),\nabla(i-1))$ is in fact lifted from a
non-zero element of $\Ext^1_A(L(i),\nabla(i-1))$ via the map, 
induced by the projection $\Delta(i)\tto L(i)$. Since 
$\nabla(i-1)$  has simple socle, it further follows 
that any non-zero element from $\Ext^1_A(L(i),\nabla(i-1))$ corresponds 
to a submodule of $I(i-1)$ with simple top $L(i)$. From this one easily 
derives that $[I(i-1):\nabla(i)]=\dim\Ext^1_A(\nabla(i),\nabla(i-1))$. 
Altogether, we obtain that $\dim\Hom_{A}(\Delta(i-1),\Delta(i))=
\dim\Ext^1_A(\nabla(i),\nabla(i-1))$. In the present section we show
that the spaces $\Hom_{A}(\Delta(i-1),\Delta(i))$ and
$\Ext^1_A(\nabla(i),\nabla(i-1))$ are connected via a non-degenerate
bilinear pairing in a natural way.

For every $i=1,\dots,n$ we fix a non-zero homomorphisms, 
$\alpha_i:\Delta(i)\to\nabla(i)$.
Remark that $\alpha_i$ is unique up to a scalar and maps the top of 
$\Delta(i)$ to the socle of $\nabla(i)$. For $j<i$ let 
$f:\Delta(j)\to\Delta(i)$ be some homomorphism and 
$\xi:\nabla(j)\overset{\beta}{\hookrightarrow} X
\overset{\gamma}{\tto}\nabla(i)$ be a short exact sequence. 
Consider the following diagram:
\begin{equation}\label{eq2.1}
\xymatrix{
0\ar[rr] && \nabla(j)\ar[rr]^{\beta} && X\ar[rr]^{\gamma} 
&& \nabla(i)\ar[rr] && 0 \\
&& \Delta(j)\ar@{-->}[u]^{\alpha_j}\ar[rrrr]^{f} && && 
\Delta(i) \ar@{-->}[u]_{\alpha_i}\ar@{=>}[ull]_{\varphi} && 
}.
\end{equation}
Since $j<i$, we have $\Ext_A^1(\Delta(i),\nabla(j))=0$ and
$\Hom_A(\Delta(i),\nabla(j))=0$. Hence 
\begin{displaymath}
\Hom_A(\Delta(i),X)\cong
\Hom_A(\Delta(i),\nabla(i)),
\end{displaymath}
which means that  $\alpha_i$ admits
a unique lifting, $\varphi:\Delta(i)\to X$, such that the triangle in
\eqref{eq2.1} commutes. Further, $L(j)$ occurs exactly once in the socle 
of $X$ and $\beta\circ \alpha_j$ is a projection of
$\Delta(j)$ onto this socle $L(j)$-component of $X$.
On the other hand, since $\gamma\circ \varphi=\alpha_i$, it follows that
$\varphi(\Rad(\Delta(i)))\subset\beta(\nabla(j))$. 
Since $[\nabla(j):L(j)]=1$, it follows that
the composition $\varphi\circ f$ is a projection of
$\Delta(j)$ onto the socle $L(j)$-component of $X$ as well.
Since $\Bbbk$ is algebraically closed, we get that
$\beta\circ \alpha_j$ and $\varphi\circ f$ differ only by a
scalar (they are not the same in general as 
$\beta\circ \alpha_j$ does depend on the choice of $\alpha_j$ and
$\varphi\circ f$ does not). Hence we can denote by  
$\langle f,\xi\rangle$ the unique element  from ${\Bbbk}$
such that $\langle f,\xi\rangle \left(\varphi\circ f\right)=
\beta\circ \alpha_j$. 

\begin{lemma}\label{l2.1}
\begin{enumerate}[(1)]
\item Let $\xi':\nabla(j)\hookrightarrow Y\tto\nabla(i)$ be a short 
exact sequence, which is congruent to $\xi$. Then 
$\langle f,\xi\rangle=\langle f,\xi'\rangle$ for any $f$ as above.
In particular, $\langle \cdot \, ,\cdot\rangle$ induces a map
from $\Hom_A(\Delta(j),\Delta(i))\times \Ext_A^1(\nabla(i),\nabla(j))$
to ${\Bbbk}$ (we will denote the induced map by the same symbol 
$\langle \cdot \, ,\cdot\rangle$ abusing notation).
\item The map $\langle \cdot \, ,\cdot\rangle:
\Hom_A(\Delta(j),\Delta(i))\times \Ext_A^1(\nabla(i),\nabla(j))\to
{\Bbbk}$ is bilinear.
\end{enumerate}
\end{lemma}

\begin{proof}
This is a standard direct calculation.
\end{proof}

Note that the form $\langle \cdot \, ,\cdot\rangle$ is independent, 
up to a non-zero scalar, of the choice of $\alpha_i$ and $\alpha_j$. 
Since the algebras $A$ and $A^{opp}$ are quasi-hereditary 
simultaneously, using the dual arguments one constructs a form,
\begin{displaymath}
\langle \cdot \, ,\cdot\rangle':
\Hom_A(\nabla(i),\nabla(j))\times \Ext_A^1(\Delta(j),\Delta(i))\to
{\Bbbk},
\end{displaymath}
and proves a dual version of  Lemma~\ref{l2.1}.

\begin{theorem}\label{t2.2}
Let $j=i-1$. Then the bilinear from $\langle \cdot \, ,\cdot\rangle$
constructed above is non-degenerate.
\end{theorem}

We remark that in the case $j<i-1$ the analogous statement is not
true in general, see the example at the beginning of Section~\ref{s3}.

\begin{proof}
First let us fix a non-zero $f:\Delta(i-1)\to\Delta(i)$. Since 
$\Delta(i-1)$ has simple top, there exists a unique submodule
$M\subset \Delta(i)$, which is maximal with respect to the condition 
$p\circ f\neq 0$,
where $p:\Delta(i)\to \Delta(i)/M$ is the natural projection. Denote 
$N=\Delta(i)/M$. The module $N$ has simple socle, which is isomorphic
to $L(i-1)$, and $p\circ f:\Delta(i-1)\to N$ is a non-zero projection
onto the socle of $N$. Now we claim that 
$\Rad(N)\hookrightarrow \nabla(i-1)$. Indeed, $\Rad(N)\subset
\Rad(\Delta(i))$ and hence it can have only composition subquotients
of the form $L(t)$, $t<i$, since $A$ is quasi-hereditary. But since
$\Rad(N)$  has the simple socle $L(i-1)$, the quasi-heredity of $A$ implies
$\Rad(N)\hookrightarrow \nabla(i-1)$ as well. Let $C$ denote the 
cokernel of this inclusion. The module $N$ is an extension of 
$\Rad(N)$ by $L(i)$ and is indecomposable. This implies that the short 
exact sequence $\xi:\Rad(N)\hookrightarrow N\tto L(i)$ represents a 
non-zero element in $\Ext_A^1(L(i),\Rad(N))$. Let us apply 
$\Hom_A(L(i),{}_-)$ to the short exact sequence 
$\Rad(N)\hookrightarrow \nabla(i-1)\tto C$
and remark that $\Hom_A(L(i),C)=0$ as $[C:L(s)]\neq 0$ implies
$s<i-1$ by above. This gives us an inclusion, 
$\Ext_A^1(L(i),\Rad(N))\hookrightarrow \Ext_A^1(L(i),\nabla(i-1))$,
and hence there exists a short exact sequence, 
$\xi':\nabla(i-1))\hookrightarrow N'\tto L(i)$, induced
by $\xi$.

\begin{lemma}\label{l2.3}
$\Ext_A^1(L(i),\nabla(i-1))\cong \Ext_A^1(\nabla(i),\nabla(i-1))$.
\end{lemma}

\begin{proof}
We apply $\Hom_A({}_-,\nabla(i-1))$ to the short exact sequence
$L(i)\hookrightarrow \nabla(i)\tto D$, where $D$ is the
cokernel of the inclusion $L(i)\hookrightarrow \nabla(i)$. 
This gives the following part in the long exact sequence: 
\begin{displaymath}
\Ext_A^{1}(D,\nabla(i-1))\to\Ext_A^1(\nabla(i),\nabla(i-1))\to
\Ext_A^1(L(i),\nabla(i-1))\to\Ext_A^{2}(D,\nabla(i-1)). 
\end{displaymath}
But $D$ contains only simple subquotients of the form
$L(s)$, $s\leq i-1$. This means that $\nabla(i-1)$ is Ext-injective
with respect to $D$ because of the quasi-heredity of $A$ and 
proves the statement.
\end{proof}

Applying Lemma~\ref{l2.3} we obtain that  the sequence $\xi'$ gives 
rise to the unique short exact sequence 
$\xi'':\nabla(i-1)\hookrightarrow N''\tto \nabla(i)$. Moreover,
by construction it also follows that $N$ is isomorphic to a submodule
in $N''$. Consider $\xi''$ with $X=N''$ in \eqref{eq2.1}. Using the
inclusion $N\hookrightarrow N''$ we obtain that the composition 
$\varphi\circ f$ is non-zero, implying $\langle f,\xi''\rangle\neq 0$. 
This proves that the left kernel of the form $\langle \cdot \, ,\cdot\rangle$ 
is zero.

To prove that  the right kernel is zero, we, basically, have to 
reverse the above arguments. Let
$\eta:\nabla(i-1)\hookrightarrow X\tto \nabla(i)$ be a non-split
short exact sequence. Quasi-heredity of $A$ implies that $\nabla(i-1)$ 
is Ext-injective with respect to $\nabla(i)/\soc(\nabla(i))$.
Hence $\eta$ is in fact a lifting of some non-split short exact
sequence, $\eta':\nabla(i-1)\hookrightarrow X'\tto L(i)$ say. In
particular, it follows that $X'$ and thus also $X$ has simple socle, 
namely $L(i-1)$. Further, applying $\Hom_A(\Delta(i),{}_-)$
to $\eta$, and using the fact that $\Delta(i)$ is Ext-projective 
with respect to $X$, one obtains that there is a unique (up to a scalar) 
non-trivial map from $\Delta(i)$ to $X$. Let $Y$ be its image. Then
$Y$ has simple top, isomorphic to  $L(i)$. Furthermore, all other simple
subquotients of $X$ are isomorphic to $L(s)$, $s<i$, and hence  $Y$ is a 
quotient of $\Delta(i)$. Since $\Delta(i-1)$ is Ext-projective with respect to
$\Rad(\Delta(i))$, we can find a map, $\Delta(i-1)\to \Rad(\Delta(i))$,
whose composition with the inclusion $\Rad(\Delta(i))\hookrightarrow
\Delta(i)$ followed by the projection from $\Delta(i)$ onto $Y$ is
non-zero. The composition of the first two maps gives us a map, 
$h:\Delta(i-1)\to\Delta(i)$, such that
$\langle h,\eta\rangle\neq 0$. Therefore the right kernel of the form 
$\langle \cdot \, ,\cdot\rangle$ is zero as well, completing the proof.
\end{proof}

\begin{corollary}\label{c2.4}
\begin{enumerate}[(1)]
\item $\Hom_A(\Delta(i-1),\Delta(i))\cong \Ext_A^1(\nabla(i),\nabla(i-1))^*$.
\item $\Ext_A^1(\Delta(i-1),\Delta(i))\cong \Hom_A(\nabla(i),\nabla(i-1))^*$.
\end{enumerate}
\end{corollary}

\begin{proof}
The first statement is an immediate corollary of Theorem~\ref{t2.2}
and the second statement follows by duality since $A^{opp}$ is 
quasi-hereditary as soon as $A$ is, see \cite{CPS1}.
\end{proof}

\begin{corollary}\label{c2.5}
Assume that $A$ has a {\em simple preserving duality}, that is a contravariant 
exact equivalence, which preserves the iso-classes of simple modules. Then 
\begin{enumerate}[(1)]
\item $\Hom_A(\Delta(i-1),\Delta(i))\cong \Ext_A^1(\Delta(i-1),\Delta(i))^*$.
\item $\Ext_A^1(\nabla(i),\nabla(i-1))\cong \Hom_A(\nabla(i),\nabla(i-1))^*$.
\end{enumerate}
\end{corollary}

\begin{proof}
Apply the simple preserving duality to the statement of Corollary~\ref{c2.4}.
\end{proof}

\section{Homomorphisms between arbitrary standard modules}\label{s3} 

It is very easy to see that the statement of Theorem~\ref{t2.2} does not
extend to the case $j<i-1$. For example, consider the path algebra $A$ of
the following quiver:
\begin{displaymath}
\xymatrix{
1 && 2\ar[ll] && 3 \ar[ll]\\
}.
\end{displaymath}
This algebra is hereditary and thus quasi-hereditary. Moreover, it is
directed and thus standard modules are projective and costandard modules are
simple. One easily obtains that $\Hom_A(\Delta(1),\Delta(3))={\Bbbk}$ 
whereas $\Ext^1_A(\nabla(3),\nabla(1))=0$. The main reason why this
happens is the fact that the non-zero homomorphism $\Delta(1)\to\Delta(3)$ 
factors through $\Delta(2)$ (note that $1<2<3$). 

Let us define another pairing in homology. Denote by $\overline{\alpha}_i$ the
natural projection of $\Delta(i)$ onto $L(i)$ and consider 
(for $j<i$) the following diagram:
\begin{equation}\label{eq3.1}
\xymatrix{
0\ar[rr] && \nabla(j)\ar[rr]^{\beta} && X\ar[rr]^{\gamma} 
&& L(i)\ar[rr] && 0 \\
&& \Delta(j)\ar@{-->}[u]^{\alpha_j}\ar[rrrr]^{f} && && 
\Delta(i) \ar@{-->}[u]_{\overline{\alpha}_i}\ar@{=>}[ull]_{\varphi} && 
}.
\end{equation}
Using this diagram, the same arguments as in Section~\ref{s2} allow us to 
define the map 
\begin{displaymath}
\overline{\langle \cdot \, ,\cdot\rangle}:
\Hom_A(\Delta(j),\Delta(i))\times \Ext_A^1(L(i),\nabla(j))\to
{\Bbbk} 
\end{displaymath}
and one can check that this map is bilinear.

\begin{proposition}\label{p3.1}
Let $N$ be the quotient of $\Delta(i)$, maximal with respect to the
following conditions: $[\Rad(N):L(s)]\neq 0$ implies $s\leq j$;
$[\Soc(N):L(s)]\neq 0$ implies $s=j$. Then
\begin{enumerate}[(1)]
\item the rank of the form
$\overline{\langle \cdot \, ,\cdot\rangle}$ equals 
the multiplicity $[\Soc(N):L(j)]$, which, in turn, is equal to
$\dim\Ext_A^1(L(i),\nabla(j))$;
\item the left kernel of $\overline{\langle \cdot \, ,\cdot\rangle}$
is the set of all morphisms $f:\Delta(j)\to\Delta(i)$ such that
$\pi\circ f=0$, where $\pi:\Delta(i)\tto N$ is the natural projection.
\end{enumerate}
\end{proposition}

\begin{proof}
The proof is analogous to that of Theorem~\ref{t2.2}.
\end{proof}

Analyzing the proof of Lemma~\ref{l2.3} it is easy to see that there is
no chance to hope for any reasonable relation between
$\Ext_A^1(L(i),\nabla(j))$ and $\Ext_A^1(\nabla(i),\nabla(j))$ in general.
However, we have the following:

\begin{proposition}\label{c3.3}
\begin{enumerate}[(1)]
\item The right kernel of $\langle \cdot \, ,\cdot\rangle$ coincides with 
the kernel of the homomorphism  $\tau:\Ext_A^1(\nabla(i),\nabla(j))\to 
\Ext_A^1(L(i),\nabla(j))$ coming from the long exact sequence in homology.
\item
Let $j=i-2$. Then $\tau$ is surjective; the rank of 
$\langle \cdot \, ,\cdot\rangle$ coincides with the rank of 
$\overline{\langle \cdot \, ,\cdot\rangle}$; and the left kernel of 
$\langle \cdot \, ,\cdot\rangle$ coincides with the left kernel of 
$\overline{\langle \cdot \, ,\cdot\rangle}$.
\end{enumerate}
\end{proposition}

\begin{proof}
The first statement follows from the proof of Theorem~\ref{t2.2}.
To prove the second statement we remark that for $j=i-2$ we have
$\Ext_A^{k}(X,\nabla(i-2))=0$, $k>1$, for any simple subquotient $X$ of
$\nabla(i)/L(i)$. This gives the surjectivity of $\tau$, which implies
all other statements.
\end{proof}

We remark that all results of this section have appropriate dual analogues.

\section{A generalization of the bilinear pairing to higher $\Ext$'s}\label{s4} 

Let us go back to the example at the beginning of Section~\ref{s3}, where
we had a hereditary algebra with $\Ext^1_A(\nabla(3),\nabla(1))=0$,
$\Hom_A(\Delta(1),\Delta(3))={\Bbbk}$, and such that any morphism from 
the last space  factors through $\Delta(2)$.
One can have the following idea: $\Hom_A(\Delta(1),\Delta(3))$ decomposes into
a product of $\Hom_A(\Delta(1),\Delta(2))$ and $\Hom_A(\Delta(2),\Delta(3))$,
by Theorem~\ref{t2.2} the space $\Hom_A(\Delta(1),\Delta(2))$ is dual to
$\Ext^1_A(\nabla(2),\nabla(1))$ and the space $\Hom_A(\Delta(2),\Delta(3))$ is 
dual to $\Ext^1_A(\nabla(3),\nabla(2))$, perhaps this means that the product of
$\Hom_A(\Delta(1),\Delta(2))$ and $\Hom_A(\Delta(2),\Delta(3))$ should correspond to 
the product of the spaces
$\Ext^1_A(\nabla(2),\nabla(1))$ and $\Ext^1_A(\nabla(3),\nabla(2))$
and thus should be perhaps paired with $\Ext^2_A(\nabla(3),\nabla(1))$ and not
$\Ext^1_A(\nabla(3),\nabla(1))$? In our example this argument does not work 
directly either since the algebra we consider is hereditary and thus $\Ext^2_A$ 
simply vanish. However, on can observe that for 
$j=i-k$, $k\in\mathbb{N}$, one could define a $\Bbbk$-linear
map from $\Ext_A^k(\nabla(i),\nabla(j))^*$ to $\Hom_A(\Delta(j),\Delta(i))$ via
\begin{multline*}
\Ext_A^k\left(\nabla(i),\nabla(j)\right)^*\overset{f}{\rightarrow}
\bigotimes_{l=0}^{k-1}
\Ext_A^1\left(\nabla(i-l),\nabla(i-l-1)\right)^*\cong 
\text{ (by Corollary~\ref{c2.4}) } \\ \cong
\bigotimes_{l=0}^{k-1}
\Hom_A^1\left(\Delta(i-l-1),\Delta(i-l)\right)\overset{g}{\rightarrow}
\Hom_A\left(\Delta(j),\Delta(i)\right),
\end{multline*}
where $g$ is the usual composition of $k$ homomorphisms, and $f$ is the dual map
to the Yoneda composition of $k$ extensions. This map would give a bilinear
pairing between $\Ext_A^k\left(\nabla(i),\nabla(j)\right)^*$ and
$\Hom_A\left(\Delta(j),\Delta(i)\right)^*$, which could also be interesting. 
However, we do not study this approach in the present paper.

Instead, we are going to  try to extend the pairing we discussed
in the previous sections to higher extensions using some resolutions. This
leads us to the following definition. Choose a minimal tilting resolution,
\begin{displaymath}
\mathcal{C}^{\bullet}:\quad\quad
0\longrightarrow T_k \overset{\varphi_k}{\longrightarrow}\dots
\overset{\varphi_2}{\longrightarrow} T_1 \overset{\varphi_1}{\longrightarrow}
T_0\overset{\varphi_0}{\longrightarrow}\nabla\longrightarrow 0,
\end{displaymath}
of $\nabla$ (see \cite[Section~5]{Ri} for the existence of such resolution).
Denote by $\mathcal{T}(\nabla)^{\bullet}$ the corresponding complex of tilting 
modules. Fix $l\in\{0,\dots,k\}$ and consider the following part of the resolution above:
\begin{displaymath}
\xymatrix{
 &&\Delta\ar@{.>}[d]^{f} &&  \\
T_{l+1}\ar[rr]^{\varphi_{l+1}} && T_l\ar@{.>}[d]^{g}\ar[rr]^{\varphi_{l}} 
&& T_{l-1} \\
 &&\nabla &&  \\
}.
\end{displaymath}
For every $f\in \Hom_A(\Delta,T_l)$ and every $g\in \Hom_A(T_l,\nabla)$
the composition gives
\begin{displaymath}
g\circ f\in \Hom_A(\Delta,\nabla)=\oplus_{i=1}^n\Hom_A(\Delta(i),\nabla(i))=
\oplus_{i=1}^n{\Bbbk}\alpha_i.
\end{displaymath}
Hence $g\circ f=\sum_{i=1}^n a_i \alpha_i$ for some $a_i\in\Bbbk$ and we can denote
$\widetilde{\langle f,g\rangle}^{(l)}= \sum_{i=1}^n a_i\in {\Bbbk}$. Obviously 
$\widetilde{\langle \cdot\, ,\cdot\rangle}^{(l)}$ defines a bilinear map from 
$\Hom_A(\Delta,T_l)\times \Hom_A(T_l,\nabla)$ to ${\Bbbk}$.
This map induces the bilinear map
\begin{displaymath}
\langle f,g\rangle^{(l)}:
\Hom_A(\Delta,T_l)\times \Hom_{Com}
\left(\mathcal{T}(\nabla)^{\bullet},\nabla^{\bullet}[l]\right)\to {\Bbbk}
\end{displaymath}
(where $\Hom_{Com}$ means the homomorphisms of complexes).
We remark that we have an obvious inclusion $\Hom_{Com}
\left(\mathcal{T}(\nabla)^{\bullet},\nabla^{\bullet}[l]\right)\subset 
\Hom_A(T_l,\nabla)$ since the complex $\nabla^{\bullet}[l]$ is 
concentrated in one degree.

\begin{theorem}\label{t4.1}
Let $f\in \Hom_A(\Delta,T_l)$ and 
$g\in \Hom_{Com}\left(\mathcal{T}(\nabla)^{\bullet},\nabla^{\bullet}[l]\right)$.
Assume that $g$ is homotopic to zero. Then $\langle f,g\rangle^{(l)}=0$. In 
particular,  $\langle \cdot\, ,\cdot\rangle^{(l)}$ induces a bilinear map,  
$\Hom_A(\Delta,T_l)\times \Ext_A^l(\nabla,\nabla)\to {\Bbbk}$. 
\end{theorem}

The form, constructed in Theorem~\ref{t4.1} will be denoted also by 
$\langle \cdot\, ,\cdot\rangle^{(l)}$ abusing notation. We remark that 
both the construction above and Theorem~\ref{t4.1} 
admit appropriate dual analogues.

\begin{proof}
Since $\mathcal{T}(\nabla)^{\bullet}$ is a complex of tilting modules, 
the second statement of the theorem follows from the first one and 
\cite[Chapter III(2), Lemma~2.1]{Ha}. To prove the first
statement we will need the following auxiliary statement.

\begin{lemma}\label{l4.2}
Let $\beta :\Delta(i)\to T(j)$ and $\gamma:T(j)\to \nabla(k)$. Then
$\gamma\circ\beta\neq 0$ if and only if $i=j=k$, $\beta\neq 0$ and 
$\gamma\neq 0$. 
\end{lemma}

\begin{proof}
Using the standard properties of tilting modules, see for example \cite{Ri}, 
we have $[T(i):L(i)]=1$, $\dim\Hom_A(\Delta(i),T(i))=1$ and any non-zero
element in this space is injective, $\dim\Hom_A(T(i),\nabla(i))=1$ and 
any non-zero element in this space is surjective. Hence in the case
$i=j=k$ the composition of non-zero $\gamma$ and $\beta$ is a non-zero
projection of the top of $\Delta(i)$ to the socle of $\nabla(i)$.
This proves the "if" statement.

To prove the "only if" statement we note that $\gamma\circ\beta\neq 0$ obviously 
implies $i=k$. Assume that $j\neq i$ and 
$\gamma\circ\beta\neq 0$.  The module $T(j)$ has a costandard filtration, 
which we fix,  and  $\Delta(i)$ is a standard module.  Hence, 
by \cite[Theorem~4]{Ri},  $\beta$ is a  linear combination of some maps, each of 
which comes from a homomorphism, which maps the top of $\Delta(i)$ to  
the socle of some $\nabla(i)$ in the costandard filtration of $T(j)$ 
(we remark that this $\nabla(i)$ is a subquotient of $T(j)$ but not a 
submodule in general).  Since the composition 
$\gamma\circ\beta$ is non-zero and $\nabla(i)$ has simple socle, we have 
that at least one whole copy of $\nabla(i)$ in the costandard filtration of 
$T(j)$ survives under $\gamma$. But, by \cite[Theorem~1]{Ri}, any costandard  
filtration of $T(j)$ ends with the subquotient $\nabla(j)\neq \nabla(i)$. 
This implies that the dimension of the image of $\gamma$ must be strictly bigger 
than $\dim \nabla(i)$, which is impossible. The obtained contradiction
shows that  $i=j=k$. The rest follows from the standard facts, used in the proof
of the "if" part.
\end{proof}

We can certainly assume that $f\in\Hom_A(\Delta(i),T_l)$ and 
$g\in\Hom_{Com}\left(\mathcal{T}(\nabla)^{\bullet},
\nabla(i)^{\bullet}[l]\right)$ for some $i$.
Consider now any homomorphism $h:T_{l-1}\to\nabla(i)$. Our aim is to show that
the composition $h\circ \varphi_l\circ f=0$.  Assume that this is not the case
and apply Lemma~\ref{l4.2} to the components of the following two 
pairs: 
\begin{enumerate}[(a)]
\item\label{ppp1} $f:\Delta(i)\to T_l$ and $h\circ \varphi_l:T_l\to \nabla(i)$
\item\label{ppp2}  $\varphi_l\circ f:\Delta(i)\to T_{l-1}$ and
$h:T_{l-1}\to \nabla(i)$.
\end{enumerate}
If $h\circ \varphi_l\circ f\neq 0$, we obtain that
both $T_l$ and $T_{l-1}$ contain a direct summand isomorphic to $T(i)$,
such that the map $\varphi_l$ induces a map, $\overline{\varphi}_l:T(i)\to T(i)$, 
which does not annihilate the unique copy of $L(i)$ inside $T(i)$. Since
$T(i)$ is indecomposable, we have that $\End_A(T(i))$ is local and thus the
non-nilpotent element $\overline{\varphi}_l\in\End_A(T(i))$ must be an
isomorphism. This contradicts the minimality of the resolution 
$\mathcal{T}(\nabla)^{\bullet}$.
\end{proof}

We remark that the sequence
\begin{displaymath}
0\to \Hom_A(\Delta,T_k)\to\dots\to \Hom_A(\Delta,T_1)\to\Hom_A(\Delta,T_0)
\to \Hom_A(\Delta,\nabla)\to 0,
\end{displaymath}
obtained from $\mathcal{C}^{\bullet}$ using $\Hom_A(\Delta,{}_-)$, is exact, and 
that Theorem~\ref{t4.1} defines a bilinear pairing between 
$\Ext_A^l(\nabla,\nabla)$ and the $l$-th element of this exact sequence. It is also 
easy to see that the pairing, given by Theorem~\ref{t4.1}, does not depend (up to 
an isomorphism of bilinear forms) on the choice of a minimal tilting resolution 
of $\nabla$. In particular, for every $l$ the rank of 
$\langle \cdot\, ,\cdot\rangle^{(l)}$ is an invariant of the algebra $A$. 
By linearity we have that 
\begin{displaymath}
\langle \cdot\, ,\cdot\rangle^{(l)}=\oplus_{i,j=1}^n
\langle \cdot\, ,\cdot\rangle^{(l)}_{i,j},
\end{displaymath}
where $\langle \cdot\, ,\cdot\rangle^{(l)}_{i,j}$ is obtained by restricting
the definition of $\langle \cdot\, ,\cdot\rangle^{(l)}$ to the homomorphisms
from $\Delta(j)$ (instead of $\Delta$) to the tilting resolution of $\nabla(i)$ 
(instead of $\nabla$). The relation between 
$\langle \cdot\, ,\cdot\rangle^{(l)}_{(i,j)}$   and the forms we have studied 
in the previous section can be described as follows: 

\begin{proposition}\label{p4.3}
$\rank \langle \cdot\, ,\cdot\rangle^{(1)}_{i,i-1}=
\dim\Ext_A^1(\nabla(i),\nabla(i-1))= \rank \langle \cdot\, ,\cdot\rangle$.
\end{proposition}

\begin{proof}
Straightforward.
\end{proof}

In the general case we have the following:

\begin{corollary}\label{c4.4}
$\rank \langle \cdot\, ,\cdot\rangle^{(l)}_{i,j}$  equals the 
multiplicity of $T(j)$ as a direct summand in the $l$-th term of 
the minimal tilting resolution of $\nabla(i)$. 
\end{corollary}

\begin{proof}
Let $T_l=\oplus_{k=1}^n T(k)^{l_k}$
and $p:T_l\tto \oplus_{k=1}^n \nabla(k)^{l_k}$ be a projection.
Since the complex $\mathcal{C}^{\bullet}$ is exact and consists of elements,
having a costandard filtration, the cokernel of any map in this complex 
has a costandard filtration itself since the category of modules with 
costandard filtration is closed with respect to taking cokernels of 
monomorphisms, see for example \cite[Theorem~1]{DR2}. This implies that
$\varphi_l$ induces a surjection from $T_l$ onto a module
having a costandard filtration. Moreover, the minimality of the 
resolution means that this surjection does not annihilate any of the
direct summands. In other words, the kernel of $\varphi_l$ is 
contained in the kernel of $p$. This implies that for the cokernel $N$ 
of $\varphi_{l+1}$ we have $\dim\Hom(N,\nabla(j))=l_j$. Using
Lemma~\ref{l4.2} it is easy to see that $\dim\Hom(N,\nabla(j))$, in fact,
equals $\rank \langle \cdot\, ,\cdot\rangle^{(l)}_{i,j}$. This 
completes the proof.
\end{proof}

We remark that, using 
Corollary~\ref{c4.4} and the Ringel duality (see \cite[Chapter~6]{Ri}), 
we can also interpret $\rank \langle \cdot\, ,\cdot\rangle^{(l)}_{i,j}$ as 
the dimension of $l$-th extension space (over $R(A)$) from the $i$-th 
standard $R(A)$-module to the $j$-th simple $R(A)$-module. For the BGG 
category $\mathcal{O}$ the dimensions of these spaces are given by 
the Kazhdan-Lusztig combinatorics.

\section{Graded non-degeneracy in a graded case}\label{s5} 

The form $\langle \cdot\, ,\cdot\rangle^{(l)}$ is degenerate in the 
general case. However, in this section we will show that it induces a  
non-degenerate pairing between the graded homomorphism and extension spaces 
for graded algebras under some assumptions in the spirit of Koszulity 
conditions. 

Throughout this section we assume that $A$ is positively graded (recall that
this means that $A=\oplus_{i\geq 0} A_i$ and $\Rad(A)=\oplus_{i> 0} A_i$).
We remark that this automatically guarantees that the simple $A$-modules 
can be considered as graded modules.  We denote by $A\mathrm{-gmod}$ the 
category of all graded (with respect to the grading fixed above) finitely 
generated $A$-modules. The morphisms in $A\mathrm{-gmod}$ are morphisms
of $A$-modules, which {\em preserve} the grading, that is these morphisms are
homogeneous morphisms of degree $0$. We denote by 
$\langle 1\rangle:A\mathrm{-gmod}\to  A\mathrm{-gmod}$ the functor, 
which shifts the grading as follows: $(M\langle 1\rangle)_i=M_{i+1}$. 

Forgetting the grading defines a faithful functor from $A\mathrm{-gmod}$
to $A\mathrm{-mod}$. We say that $M\in A\mathrm{-mod}$ admits the {\em 
graded lift} $\tilde{M}\in A\mathrm{-gmod}$ (or, simply, is {\em gradable}) 
provided that, after forgetting the grading, the 
module $\tilde{M}$ becomes isomorphic to $M$. If $M$ is indecomposable and 
admits a graded lift, then this lift is unique up to an isomorphism in 
$A\mathrm{-gmod}$  and a shift of grading, see for example 
\cite[Lemma~2.5.3]{BGS}.

For $M,N\in A\mathrm{-gmod}$ we set $\ext_A^{i}(M,N)=
\Ext_{A\mathrm{-gmod}}^{i}(M,N)$, $i\geq 0$. It is clear that, forgetting the 
grading, we have 
\begin{equation}\label{greq}
\Ext_A^{i}(M,N)=\oplus_{j\in\Z}\ext_A^{i}(M,N\langle j\rangle), \quad\quad
i\geq 0
\end{equation}
(see for example \cite[Lemma~3.9.2]{BGS}). 

\begin{lemma}\label{l5.1}
Let $M,N\in A\mathrm{-gmod}$. Then the non-graded trace $\mathrm{Tr}_M(N)$ of 
$M$ in $N$, that is the sum of the images of all (non-graded) homomorphism 
$f:M\to N$, belongs to $A\mathrm{-gmod}$.
\end{lemma}

\begin{proof}
Any $f:M\to N$ can be written as a sum of homogeneous components 
$f_i:M\to N\langle i\rangle$, $i\in\Z$, in particular, the image of $f$ is 
contained in the sum of the images of all $f_i$. Since the image of a 
homogeneous map is a graded submodule of $N$, the statement follows.
\end{proof}

\begin{corollary}\label{c5.2}
All  standard and costandard  $A$-modules are gradable.
\end{corollary}

\begin{proof}
By duality it is enough to prove the statement for standard modules.
The module $\Delta(i)$ is defined as a quotient of $P(i)$ modulo the
trace of $P(i+1)\oplus\dots\oplus P(n)$ in $P(0)$. For positively graded
algebras all projective modules are obviously graded and hence the 
statement follows from Lemma~\ref{l5.1}.
\end{proof}

\begin{proposition}\label{p5.3}
Let $M,N\in A\mathrm{-gmod}$. Then the universal extension of $M$ by
$N$ (in the category $A\mathrm{-mod}$) is gradable.
\end{proposition}

\begin{proof}
As we have mentioned before, we have 
$\Ext_A^{1}(M,N)=\oplus_{j\in\Z}\ext_A^{1}(M,N\langle j\rangle)$. 
Every homogeneous extension obviously produces a gradable module. 
Since we can construct the universal extension
of $N$ by $M$ choosing a homogeneous basis in $\Ext_A^1(M,N)$, the 
previous argument shows that the obtained module will be gradable. 
This completes the proof.
\end{proof}

We would like to fix a grading on all modules, related to 
the quasi-hereditary structure. We concentrate $L$ in degree 
$0$ and fix the gradings on $P$, $\Delta$, $\nabla$ and $I$ such that the
canonical maps $P\tto L$, $\Delta\tto L$, $L\hookrightarrow \nabla$
and $L\hookrightarrow I$ are all morphism in $A\mathrm{-gmod}$. The
only structural modules, which are left, are tilting modules.
However, to proceed, we have to show first that tilting modules are gradable. 

\begin{corollary}\label{c5.4}
All tilting $A$-modules admit graded lifts. Moreover, for $T$ this lift can be
chosen such that both the inclusion $\Delta\hookrightarrow T$ and the 
projection $T\tto \nabla$ are morphisms in $A\mathrm{-gmod}$.
\end{corollary}

\begin{proof}
By \cite[Proof of Lemma~3]{Ri}, the tilting $A$-module $T(i)$  is produced 
by a sequence of universal extensions as follows:  we start from the 
(gradable) module $\Delta(i)$, and on each step we extend some (gradable) 
module $\Delta(j)$, $j<i$, with the module, obtained on the previous step. 
Using  Proposition~\ref{p5.3} and induction we see that all modules, obtained 
during this process, are gradable. The statement about the choice of the
lift is obvious.
\end{proof}

We fix the grading on $T$, given by Corollary~\ref{c5.4}. This automatically
induces a grading on the Ringel dual $R(A)=\End_A(T)^{opp}$. In what follows 
we always will consider $R(A)$ as a graded algebra with respect to this 
induced  grading.

Note that the 
same ungraded $A$-module can occur as a part of different structures, for 
example, a module can be projective, injective and tilting at the same time. 
In this case it is possible that the lifts of this module, which we fix for
different structures, are different. For example, if we have a non-simple
projective-injective module, then, considered as a projective module, it is
graded in non-negative degrees with top being in degree $0$; considered as 
an injective module, it is graded in non-positive degrees with socle being 
in degree $0$; and, considered as a tilting module, it has non-trivial 
components both in negative and positive degrees.

A complex, $\mathcal{X}^{\bullet}$, of graded projective (resp. 
injective, resp. tilting) modules will be called {\em linear} provided 
that $\mathcal{X}^{i}\in\mathrm{add} (P\langle i\rangle)$
(resp. $\mathcal{X}^{i}\in\mathrm{add} (I\langle i\rangle)$, resp.
$\mathcal{X}^{i}\in\mathrm{add} (T\langle i\rangle)$) for all $i\in\Z$.

To avoid confusions between the degree of a graded component of a module and 
the degree of a component in some complex, to indicate the place of a component 
in a complex we will use the word {\em position} instead of the word degree.

We say that $M\in A\mathrm{-gmod}$ admits an {\em LT-resolution}, 
$\mathcal{T}^{\bullet}\tto M$, (here LT stands for linear-tilting) if 
$\mathcal{T}^{\bullet}$ is a linear complex of tilting modules from 
$A\mathrm{-gmod}$, such that $\mathcal{T}^{i}=0$, $i>0$, and the 
homology of $\mathcal{T}^{\bullet}$ is concentrated in position $0$ 
and equals $M$ in this position. One also defines {\em LT-coresolution} 
in the dual way. The main result of this section  is the following:

\begin{theorem}\label{t5.5}
Let $A$ be a positively graded quasi-hereditary algebra and $1\leq i,j\leq n$. 
Assume that 
\begin{enumerate}[(i)]
\item\label{l5.5.1} $\nabla(i)$ admits an LT-resolution, 
$\mathcal{T}(\nabla(i))^{\bullet}\tto\nabla(i)$;
\item\label{l5.5.2} the induced grading on  $R(A)$ is positive. 
\end{enumerate}
Then the form $\langle \cdot\, ,\cdot\rangle^{(l)}_{i,j}$ induces a 
non-degenerate bilinear pairing between 
\begin{displaymath}
\hom_A(\Delta(j)\langle -l\rangle, \mathcal{T}(\nabla(i))^{-l})\quad 
\text{ and }\quad  \ext_A^l(\nabla(i),\nabla(j)\langle -l\rangle).
\end{displaymath}
\end{theorem}

We remark that Theorem~\ref{t5.5} has a dual analogue.

\begin{proof}
The assumption \eqref{l5.5.2} means that 
\begin{gather}
\hom_A(\Delta\langle s\rangle,T)\neq 0 \quad\quad\mathrm{ implies }
\quad\quad s\leq 0 \label{bl1}\\
\hom_A(\Delta(k)\langle s\rangle,T(m))\neq 0 \quad\mathrm{ and }\quad k\neq m
\quad\quad\mathrm{ implies }\quad\quad s<0\label{bl2}.
\end{gather}
Hence, it follows that 
\begin{displaymath}
\dim \hom_A\left(\Delta(j)\langle -l\rangle, \mathcal{T}(\nabla(i))^{-l}\right)
\end{displaymath}
equals the multiplicity of $T(j)\langle -l\rangle$ as a direct summand of
$\mathcal{T}(\nabla(i))^{-l}$, which, using the dual arguments, in turn, equals 
\begin{displaymath}
\dim \hom_A\left(\mathcal{T}(\nabla(i))^{-l},\nabla(j)\langle -l\rangle\right). 
\end{displaymath}
From the definition of an LT-resolution and \eqref{bl1}-\eqref{bl2} we also obtain 
\begin{displaymath}
\hom_A\left(\mathcal{T}(\nabla(i))^{-l+1},\nabla(j)\langle -l\rangle\right)=0,
\end{displaymath}
which means that there is no homotopy from $\mathcal{T}(\nabla(i))^{\bullet}$
to $\nabla(j)\langle -l\rangle^{\bullet}$. The arguments, analogous to those, 
used in Corollary~\ref{c4.4}, imply that any map from 
$\mathcal{T}(\nabla(i))^{-l}$ to $\nabla(j)\langle -l\rangle$ induces a
morphism of complexes from $\mathcal{T}(\nabla(i))^{\bullet}$ to 
$\nabla(j)\langle -l\rangle^{\bullet}$. Hence 
\begin{displaymath}
\dim \ext_A^l\left(\nabla(i),\nabla(j)\langle -l\rangle\right)= 
\dim \hom_A\left(\mathcal{T}(\nabla(i))^{-l},\nabla(j)\langle -l\rangle\right).
\end{displaymath}

We can now interpret every $f\in \hom_A\left(\Delta(j)\langle -l\rangle, 
\mathcal{T}(\nabla(i))^{-l}\right)$ as a fixation of a direct summand of 
$\mathcal{T}(\nabla)^{-l}$, which is isomorphic to $T(i)\langle -l\rangle$. 
Projecting it further onto  $\nabla(j)\langle -l\rangle$ shows that 
the left kernel of the form $\langle \cdot\, ,\cdot\rangle^{(l)}_{i,j}$
is zero. Since the dimensions of the left and the right spaces coincide 
by the arguments above, we obtain that the form is
non-degenerate. This completes the proof.
\end{proof}

It is easy to see that the condition \eqref{l5.5.2} of Theorem~\ref{t5.5} does
not imply the condition \eqref{l5.5.1} in general. Further, it is also easy to
see, for example for the path algebra of the following quiver:
\begin{displaymath}
\xymatrix{ 1 && 2\ar[rr] && 3\ar@/_1pc/[llll] && 4\ar[ll]},
\end{displaymath}
that the condition \eqref{l5.5.1} (even if we assume it to be satisfied for 
all $i$) does not imply the condition \eqref{l5.5.2} in general, 
However, we do not know if the assumptions of the existence of an
$LT$-resolution for $\nabla$ and, simultaneously, an $LT$-coresolution
for $\Delta$, would imply the condition \eqref{l5.5.2}.
 
We also would like to remark that the conditions of Theorem~\ref{t5.5} are 
not at all automatic even in very good cases. For example one can check that the 
path algebra of the following quiver:
\begin{displaymath}
\xymatrix{ 1 && 2\ar[ll]\ar[rr] && 3\ar[rr] && 4}
\end{displaymath}
is standard Koszul, however, both conditions of Theorem~\ref{t5.5} fail. 

Let $A$ be a positively graded quasi-hereditary algebra. We say that
$A$ is an {\em SCK-algebra} (abbreviating standard-costandard-Koszul) 
provided that $A$ is standard Koszul and the induced grading on $R(A)$
is positive. We say that $A$ is an {\em SCT-algebra} (abbreviating 
standard-costandard-tilting) provided that all standard and costandard
modules admit LT-(co)resolutions. By \cite[theorem~1]{ADL}, any standard 
Koszul algebra, and thus any SCK-algebra, is Koszul. We finish this 
section with the following observation.

\begin{theorem}\label{t5.6}
Any SCK-algebra is an SCT-algebra and vice versa.
\end{theorem}

\begin{proof}
Our first observation is that for any SCT-algebra $A$ the induced grading
on the $R(A)$ is positive. To prove this it is enough to show that all
subquotients in any standard filtration of the cokernel of the morphism
$\Delta(i)\hookrightarrow T(i)$ have the form $\Delta(j)\langle l\rangle$, 
$l>0$. This follows by induction in $i$. For $i=1$ the statement is obvious,
and the induction step follows from the inductive assumption applied
to the first term of the linear tilting coresolution of $\Delta(i)$.

Now we claim that the Ringel dual of an SCT algebra is SCK and vice versa.
Assume that $A$ is SCT. Applying $\Hom_A(T,{}_-)$ to the LT-resolution
of $\nabla$ we obtain that the $k$-th component of the projective 
resolution of the standard  $R(A)$-module is generated in degree $k$.
Applying analogous arguments to the LT-coresolution of $\Delta$ we obtain 
that the $k$-th component of the injective resolution of the costandard  
$R(A)$-module is generated in degree $-k$. As we have already shown, the
induced grading on $R(A)$ is positive. Furthermore, the (graded) Ringel 
duality maps injective $A$-modules to tilting $R(A)$-modules, which implies
that the grading, induced on $A$ from $R(A)\mathrm{-gmod}$, will 
coincide with the original grading on $A$, and hence will be positive 
as well. This means that $R$ is SCK. The arguments in the opposite direction 
are similar. 

To complete the proof it is now enough to show, say, that any SCT algebra is
SCK. The existence of a linear tilting coresolution for $\Delta$ and the
above proved fact that for an SCT-algebra $A$ the induced grading
on the $R(A)$ is positive, imply $\ext^{k}(\Delta\langle l\rangle,\Delta)=0$
unless $l\leq k$. Since  $A$ is positively graded, we have that the 
$k$-th term of the projective resolution of $\Delta$ consists of 
modules of the form $P(i)\langle -l\rangle$, $l\geq k$. Assume that
for some $k$ we have that $P(i)\langle -l\rangle$ with $l>k$ occurs.
Since every kernel and cokernel in our resolution has a standard
filtration, we obtain that $\ext^{k}(\Delta,\Delta(i)\langle -l\rangle)\neq 0$
with $l> k$, which contradicts $l\leq k$ above. This implies that 
$\Delta$ has a linear projective resolution. Analogous arguments imply that
$\nabla$ has a linear injective coresolution. This completes the proof.
\end{proof}

\section{The category of linear complexes of tilting modules}\label{s55} 

We continue to work under the assumptions of Section~\ref{s5}, moreover, 
we assume, until the end of this section, that $A$ is such that both 
$A$ and $R(A)$ are positively graded.

The results of Section~\ref{s5} motivate the following definition:
We say that $M\in A\mathrm{-gmod}$ is {\em $T$-Koszul}
provided that $M$ is isomorphic in $D^b(A\mathrm{-gmod})$ to a linear 
complex of tilting modules. Thus any module, which admits an 
$LT$-(co)resolution, is $T$-Koszul.

We denote by $\mathfrak{T}=\mathfrak{T}(A)$ the category, whose objects are
linear complexes of tilting modules and morphisms are all
morphisms of graded complexes (which means that all components of these
morphisms are homogeneous homomorphisms of $A$-modules of degree $0$). 
We also denote by $\mathfrak{T}^b$ the full subcategory of 
$\mathfrak{T}$, which consists of bounded complexes.

\begin{lemma}\label{l55.1}
\begin{enumerate}[(1)]
\item $\mathfrak{T}$ is an abelian category.
\item $\langle -1\rangle[1]:\mathfrak{T}\to \mathfrak{T}$ is 
an auto-equivalence.
\item The complexes $\left(T(i)^{\bullet}\right)\langle -l\rangle[l]$
constitute an exhaustive list of simple objects in $\mathfrak{T}$.
\end{enumerate}
\end{lemma}

\begin{proof}
The assumption that the grading on $R(A)$, induced from $A\mathrm{-gmod}$,
is positive, implies that the algebra $\mathrm{end}_A(T^{\bullet})$ is 
semi-simple. Using this it is easy to check that taking the usual kernels 
and cokernels of morphisms of complexes defines on $\mathfrak{T}$ the 
structure of an abelian category. That 
$\langle -1\rangle[1]:\mathfrak{T}\to \mathfrak{T}$ is an auto-equivalence
follows from the definition.

The fact that $\mathrm{end}_A(T^{\bullet})$ is semi-simple and the above 
definition of the abelian structure on $\mathfrak{T}$ imply that any 
non-zero homomorphism in $\mathfrak{T}$ to the complex 
$\left(T(i)^{\bullet}\right)\langle -l\rangle[l]$  is surjective. Hence 
the objects $\left(T(i)^{\bullet}\right)\langle -l\rangle[l]$ are simple. 
On the other hand, 
it is easy to see that for any linear complex $\mathcal{T}^{\bullet}$
and for any $k\in\Z$ the complex $\left(\mathcal{T}^{k}\right)^{\bullet}$
is a subquotient of $\mathcal{T}^{\bullet}$ provided that 
$\mathcal{T}^{k}\neq 0$. Hence any simple object in $\mathfrak{T}$ should 
contain only one non-zero component. In order to be a simple object, this 
component obviously should be an indecomposable $A$-module. Therefore any
simple object in $\mathfrak{T}$ is isomorphic to 
$\left(T(i)^{\bullet}\right)\langle -l\rangle[l]$ for some $i$ and $l$. 
This  completes the proof.
\end{proof}

Our aim is to show that $\mathfrak{T}$ has enough projective objects.
However, to do this it is more convenient to switch to a different 
language and to prove a more general result. 

Let $B=\oplus_{i\in\Z}B_i$ be a basic positively graded $\Bbbk$-algebra such 
that $\dim_{\Bbbk}B_i<\infty$ for all $i\geq 0$. Denote by $\mathfrak{B}$ 
the category of linear complexes of projective $B$-modules, and by 
$\tilde{\mathfrak{B}}$ the category, whose objects are all sequences 
$\mathcal{P}^{\bullet}$ of projective $B$-modules, such that 
$\mathcal{P}^{i}\in\mathrm{add}(P\langle -i\rangle)$ for all $i\in\Z$,
and whose morphisms are all morphisms of graded sequences (consisting 
of homogeneous maps of degree $0$). The objects of
$\tilde{\mathfrak{B}}$ will be called {\em linear sequences of projective 
modules}.

Denote by $\mu: B_1\otimes_{B_0}B_1\to B_2$ the multiplication 
map and by $\mu^{*}:B_2^*\to B_1^*\otimes_{B_0^*}B_1^*$ the dual map.
Define the algebra $\Lambda$ as the quotient of the free
positively graded tensor algebra $B_0[B_1^*]$ modulo the homogeneous 
ideal generated by $\mu^*(B_2^*)$. 

A graded module, $M=\oplus_{i\in\Z}M_i$, over a graded algebra is 
called {\em locally finite} provided that $\dim M_i<\infty$ for all $i$. 
Note that a locally finite module does not need to be finitely generated.
For a graded algebra, $C$, we denote by $C\mathrm{-lfmod}$ the category 
of all locally finite graded $C$-modules (with morphisms being homogeneous
maps of degree $0$).

The following statement was proved in \cite[Theorem~2.4]{MS}. For the
sake of completeness we present a short version of the proof.

\begin{theorem}\label{t55.2}
There is an equivalence of categories, $\overline{F}: 
B_0[B_1^*]\mathrm{-lfmod}\to \tilde{\mathfrak{B}}$,  
which induces an equivalence,
$F:\Lambda\mathrm{-lfmod}\to \mathfrak{B}$.
\end{theorem}

\begin{proof}
Let $P$ denote the projective generator of $B$.
We construct the functor $\overline{F}$ in the following way: 
Let $X=\oplus_{j\in\Z}X_j\in B_0[B_1^*]\mathrm{-lfmod}$. We define 
$\overline{F}(X)=\mathcal{P}^{\bullet}$, where
$\mathcal{P}^{j}=P\langle j\rangle \otimes_{B_0} X_{j}$, $j\in \Z$.
To define the differential $d_j:\mathcal{P}^{j}\to \mathcal{P}^{j+1}$
we note that $P\cong {}_B B$ and use the following bijections:
\begin{equation}\label{eq55.4}
\begin{array}{lcl}
\displaystyle
\{M\in B_0[B_1^*]\mathrm{-lfmod}: 
M|_{B_0}= X|_{B_0} \} & \cong &
\text{ (since $B_0[B_1^*]$ is free) } \\ \displaystyle
\prod_{j\in\Z}\mathrm{hom}_{B_0-B_0}\left(B_1^*\langle j+1\rangle,
\mathrm{Hom}_{\Bbbk}\left(X_{j},X_{j+1}\right)\right)
& \cong & \text{ (by adjoint associativity)} \\\displaystyle
\prod_{j\in\Z}\mathrm{hom}_{B_0}\left(X_{j},
B_1\langle j+1\rangle\otimes_{B_0} X_{j+1}, \right)
& \cong & \text{ (because of grading) } \\ \displaystyle
\prod_{j\in\Z}\mathrm{hom}_{B_0}\left(X_{j},
B\langle j+1\rangle\otimes_{B_0} X_{j+1}, \right) &\cong &
\text{ (by projectivity of ${}_B B$)}\\ \displaystyle
\prod_{j\in\Z}\mathrm{hom}_B\left(B\langle j\rangle\otimes_{B_0} X_{j},
B\langle j+1\rangle\otimes_{B_0} X_{j+1} \right). & & \\
\end{array}
\end{equation}
Thus, starting from the fixed $X$, the equalities of \eqref{eq55.4} 
produce for each $j\in Z$ a unique map from the space  
$\mathrm{hom}_B\left(B\langle j\rangle\otimes_{B_0} X_{j},
B\langle j+1\rangle\otimes_{B_0} X_{j+1} \right)$, which defines the 
differential in $\mathcal{P}^{\bullet}$.

Tensoring with the identity map on ${}_B B$ the correspondence
$\overline{F}$, defined above on objects, extends to a functor 
from $B_0[B_1^*]\mathrm{-lfmod}$ to $\tilde{\mathfrak{B}}$. Since 
$\hom({}_B B,{}_B B)\cong B_0$ is a direct sum of several copies
of $\Bbbk$, it follows by a direct calculation that 
$\overline{F}$ is full and faithful. It is also easy to derive from the
construction that $\overline{F}$ is dense. Hence it is an equivalence
of categories $B_0[B_1^*]\mathrm{-lfmod}$ and $\tilde{\mathfrak{B}}$.

Now the principal question is: when $\overline{F}(X)$ is a complex? Let
\begin{gather*}
d_j: B\langle j\rangle\otimes_{B_0} X_{j}\to 
B\langle j+1\rangle\otimes_{B_0} X_{j+1},\\
d_{j-1}: B\langle j-1\rangle\otimes_{B_0} X_{j-1}\to 
B\langle j\rangle\otimes_{B_0} X_{j}
\end{gather*}
be as constructed above.
Let further
\begin{gather*}
\delta_{j}:X_{j}\to B_1 \langle j+1\rangle\otimes_{B_0} X_{j+1},\\
\delta_{j-1}:X_{j-1}\to B_1 \langle j\rangle\otimes_{B_0} X_{j}
\end{gather*}
be the corresponding maps, given by \eqref{eq55.4}. Then
$d_jd_{j-1}=0$ if and only if
\begin{displaymath}
\left(\mu\otimes \mathrm{Id}_{X_{j+1}}\right)\circ 
\left(\mathrm{Id}_{B_1}\otimes\delta_j\right)\circ\delta_{j-1}=0.
\end{displaymath}
The last equality, in turn, is equivalent to the fact that the
global composition of morphisms in the 
following diagram is zero:
\begin{displaymath}
B_2^* \xrightarrow{\mu^*}
B_1^*\otimes B_1^* \xrightarrow{b}
\mathrm{Hom}_{\Bbbk}\left(X_{j},X_{j+1}\right)\otimes
\mathrm{Hom}_{\Bbbk}\left(X_{j-1},X_{j}\right)\xrightarrow{c}
\mathrm{Hom}_{\Bbbk}\left(X_{j-1},X_{j+1}\right),
\end{displaymath}
where the map $b$ is given by two different applications of 
\eqref{eq55.4} and $c$ denotes the 
usual composition. Hence $\overline{F}(X)$ is a complex
if and only if $\mathrm{Im}(\mu^*) X=0$ or, equivalently,
$X\in \Lambda\mathrm{-lfmod}$.
\end{proof}

It is clear that the equivalence, constructed in the proof of 
Theorem~\ref{t55.2}, sends the auto-equivalence $\langle 1\rangle$ on
$\Lambda\mathrm{-lfmod}$ (resp. on $B_0[B_1^*]\mathrm{-lfmod}$) to 
the auto-equivalence $\langle -1\rangle[1]$
on $\mathfrak{B}$ (resp. on $\tilde{\mathfrak{B}}$).

Now we are back to the original setup of this section.

\begin{corollary}\label{c55.6}
Let $R=R(A)=\oplus_{i\geq 0}R_i$ and set 
$\Lambda=R_0[R_1^*]/(\mu^*(R_2^*))$,
where $\mu$ denotes the multiplication in $R$. Then the category 
$\mathfrak{T}$ is equivalent to $\Lambda\mathrm{-lfmod}$.
\end{corollary}

\begin{proof}
Apply first the graded Ringel duality and then Theorem~\ref{t55.2}.
\end{proof}

\begin{corollary}\label{c55.7}
Assume that $R=R(A)$ is Koszul. Set 
$\Lambda=(E(R(A)))^{opp}$. Then
\begin{enumerate}[(1)]
\item $\mathfrak{T}$ is equivalent to the category $\Lambda\mathrm{-lfmod}$. 
\item The category $\mathfrak{T}^b$ is equivalent to $\Lambda\mathrm{-gmod}$.
\end{enumerate}
\end{corollary}

\begin{proof}
If the algebra $R=\oplus_{i\geq 0}R_i$ is Koszul then, 
by \cite[Section~2.9]{BGS}, the 
formal quadratic dual algebra $R_0[R_1^*]/(\mu^*(R_2^*))$ is isomorphic
to $(E(R))^{opp}$. Now everything  follows from  Corollary~\ref{c55.6}.
\end{proof}

Corollary~\ref{c55.6} motivates the further study of the categories
$\mathfrak{T}$ and $\mathfrak{T}^b$. We start with a description of 
the first extension spaces between the simple objects in $\mathfrak{T}$. 
Surprisingly enough, this result can be obtained without any 
additional assumptions.

\begin{lemma}\label{l55.8}
Let $i,j\in\{1,\dots,n\}$ and $l\in \Z$. Then
$\mathrm{ext}_{\mathfrak{T}}^1\left(
T(i)^{\bullet},T(j)^{\bullet}\langle -l\rangle[l]\right)\neq 0$
implies $l=-1$. Moreover, 
\begin{displaymath}
\mathrm{ext}_{\mathfrak{T}}^1\left(
T(i)^{\bullet},T(j)^{\bullet}\langle 1\rangle[-1]\right)\cong
\mathrm{hom}_A\left(T(i),T(j)\langle 1\rangle\right).
\end{displaymath}
\end{lemma}

\begin{proof}
A direct calculation, using the definition of the first extension
via short exact sequences and the abelian structure on $\mathfrak{T}$.
\end{proof}

Recall from \cite{CPS1,DR} that an associative algebra is quasi-hereditary 
if and only if its module category is a highest weight category.
Our goal is to establish some conditions under which $\mathfrak{T}^b$
becomes a highest weight category. To prove that a category is a 
highest weight category one has to determine the (co)standard objects.

\begin{proposition}\label{p55.9}
\begin{enumerate}[(1)]
\item Assume that $\Delta(i)$ admits an LT-coresolution, 
$\Delta(i)\hookrightarrow \mathcal{T}(\Delta(i))^{\bullet}$, for all $i$. 
Then $\mathrm{ext}_{\mathfrak{T}}^1\left(\mathcal{T}(\Delta(i))^{\bullet},
T(j)^{\bullet}\langle -l\rangle[l]\right)= 0$
for all $l\in\Z$ and $j\leq i$.
\item Assume that $\nabla(i)$ admits an LT-resolution, 
$\mathcal{T}(\nabla(i))^{\bullet}\tto \nabla(i)$, for all $i$. Then we have
$\mathrm{ext}_{\mathfrak{T}}^1\left(T(j)^{\bullet}\langle -l\rangle[l],
\mathcal{T}(\nabla(i))^{\bullet}\right)= 0$
for all $l\in\Z$ and $j\leq i$.
\end{enumerate}
\end{proposition}

\begin{proof}
By duality, it is certainly enough to prove only the first statement.
Using the induction with respect to the quasi-hereditary structure it
is even enough to show that $\mathcal{T}(\Delta(n))^{\bullet}$ is
projective in $\mathfrak{T}$. By Lemma~\ref{l55.8} we can also assume 
that $l<0$. Let 
\begin{equation}\label{eq55.9.1}
0\to T(j)^{\bullet}\langle -l\rangle[l]\to \mathcal{X}^{\bullet}
\to \mathcal{T}(\Delta(n))^{\bullet}\to 0
\end{equation}
be a short exact sequence in $\mathfrak{T}$. Let further  $d^{\bullet}$
denote the differential in $\mathcal{T}(\Delta(n))^{\bullet}$. Consider
the short exact sequence
\begin{equation}\label{eq55.9.2}
0\to \ker(d^{-l})\to \mathcal{T}(\Delta(n))^{-l}\to \ker(d^{-l+1})\to 0.
\end{equation}
Since $\mathcal{T}(\Delta(n))^{\bullet}$ is a tilting coresolution of 
a standard module, it follows that all modules in \eqref{eq55.9.2} have 
standard filtration. Hence, applying $\Hom_A({}_-,T(j))$ to 
\eqref{eq55.9.2}, and using the fact that
$T(j)$ has a costandard filtration, we obtain the surjection
\begin{displaymath}
\Hom_A\left(\mathcal{T}(\Delta(n))^{-l},T(j)\right)\tto
\Hom_A\left(\ker(d^{-l}),T(j)\right),
\end{displaymath}
which induces the graded surjection
\begin{displaymath}
\hom_A\left(\mathcal{T}(\Delta(n))^{-l},T(j)\langle -l\rangle\right)\tto
\hom_A\left(\ker(d^{-l}),T(j)\langle -l\rangle\right).
\end{displaymath}
The last surjection allows one to perform a base change in 
$\mathcal{X}^{\bullet}$, which splits the sequence \eqref{eq55.9.1}.
This proves the statement.
\end{proof}

For $R=R(A)$ we
introduce the notation $R^{!}=R_0[R_1^*]/(\mu^*(R_2^*))$
(if $R$ is Koszul, this notation coincides 
with the one used for the formal quadratic dual in \cite[2.8]{BGS}). 
We have $\mathfrak{T}\cong R^{!}\mathrm{-lfmod}$ by Corollary~\ref{c55.6}. 

\cite[Theorem~1]{ADL} states that a standard Koszul quasi-hereditary
algebra is Koszul (which means that if standard modules admit linear projective
resolutions and costandard modules admit linear injective resolutions, then
simple modules admit both linear projective and linear injective resolutions). 
An analogue of this statement in our case is the  following:
 
\begin{theorem}\label{t55.12}
Let $A$ be an SCT algebra. Then
\begin{enumerate}[(1)]
\item\label{t55.12.1} $R^{!}$ is quasi-hereditary with respect to the
usual order on $\{1,2,\dots,n\}$, or, equivalently,
$\mathfrak{T}^b\simeq R^{!}\mathrm{-gmod}$ is a highest weight category;
\item\label{t55.12.2} $\mathcal{T}(\Delta(i))^{\bullet}$, $i=1,\dots,n$, are standard
objects in $\mathfrak{T}^b$;
\item\label{t55.12.3}$\mathcal{T}(\nabla(i))^{\bullet}$, $i=1,\dots,n$, are costandard
objects in $\mathfrak{T}^b$;
\end{enumerate}
Assume further that the algebra $R^{!}$ is SCT. Then
\begin{enumerate}[(1)]
\setcounter{enumi}{3}
\item\label{t55.12.4} simple $A$-modules are $T$-Koszul, in particular, for every
$i=1,\dots,n$ there exists a linear complex, $\mathcal{T}(L(i))^{\bullet}$, 
of tilting modules, which is isomorphic to $L(i)$ in $D^b(A\mathrm{-gmod})$;
\item\label{t55.12.5} $\mathcal{T}(L(i))^{\bullet}$, $i=1,\dots,n$, are tilting
objects with respect to the quasi-hereditary structure on $\mathfrak{T}^b$.
\end{enumerate}
\end{theorem}

\begin{proof}
The algebra $R$ is quasi-hereditary with respect to the opposite order
on $\{1,2,\dots,n\}$. Moreover, $R$ is SCK by Theorem~\ref{t5.6}, in particular,
it is standard Koszul, thus also Koszul by \cite[Theorem~1]{ADL}. Hence its 
Koszul dual, which is isomorphic to $(R^{!})^{opp}$ by \cite[2.10]{BGS}, 
is quasi-hereditary with respect to the usual order on $\{1,2,\dots,n\}$
by \cite[Theorem~2]{ADL}. This certainly means that $R^{!}$ is 
quasi-hereditary with respect to the usual order on $\{1,2,\dots,n\}$.
From Corollary~\ref{c55.7} we also obtain $R^{!}\mathrm{-gmod}\simeq
\mathfrak{T}^b$. This proves the first statement. 

That the objects $\mathcal{T}(\Delta(i))^{\bullet}$, $i=1,\dots,n$, are 
standard and the objects $\mathcal{T}(\nabla(i))^{\bullet}$, $i=1,\dots,n$, 
are costandard follows from Proposition~\ref{p55.9} and \cite[Theorem~1]{DR2}.
This proves \eqref{t55.12.2} and \eqref{t55.12.3}.

Now we can assume that $R^{!}$ is an SCT-algebra. In particular, it is
quasi-hereditary, and hence the category 
$\mathfrak{T}^b$ must contain tilting objects with respect to the
corresponding highest weight structure. By \cite[Proof of Lemma~3]{Ri}, 
the tilting objects in $\mathfrak{T}^b$ can be constructed via a 
sequence of universal extensions, which starts with some standard 
object and proceeds by extending other (shifted) standard objects 
by objects, already constructed on previous steps. The assumption
that $R^{!}$ is SCT=SCK means that new standard objects should be shifted by 
$\langle -l\rangle[l]$ with $l>0$. From the second statement 
of our theorem, which we have already proved above, it follows that the 
standard objects in $\mathfrak{T}^b$ are exhausted by 
$\mathcal{T}(\Delta(i))^{\bullet}$,  $i=1,\dots,n$, and their shifts. 
The homology of $\mathcal{T}(\Delta(i))^{\bullet}$ is concentrated in 
position $0$ and in non-negative degrees. It follows that the homology 
of the tilting object in $\mathfrak{T}^b$, which we obtain, using this
construction, will be concentrated in non-positive positions and in 
non-negative degrees. 

On the other hand, a dual construction, that is the one, which uses 
costandard objects, implies that the homology of the same tilting 
object in $\mathfrak{T}^b$
will be concentrated in non-negative positions and in non-positive degrees.
This means that the homology of an indecomposable  tilting object in 
$\mathfrak{T}^b$ is concentrated in position $0$ and in degree $0$ and
hence is a simple $A$-module. This proves two last statements of our
theorem and completes the proof.
\end{proof}

In the next section we will show that all the above conditions are satisfied
for the associative algebras, associated with the blocks of the 
BGG category $\mathcal{O}$.

We remark that, under conditions of Theorem~\ref{t55.12}, in the category 
$\mathfrak{T}$ the standard and costandard $A$-modules remain standard and 
costandard objects respectively via their tilting (co)resolutions. Tilting
$A$-modules become simple objects, and simple $A$-modules become tilting
objects via $\mathcal{T}(L(i))^{\bullet}$.  

An SCT algebra $A$ for which $R(A)^!$ is SCT will be called {\em balanced}.
The results of this section allow us to formulae a new type of duality 
for balanced algebras (in fact, this just means that we can perform in one 
step the following path $A\leadsto R\leadsto R^{!}\leadsto R(R^{!})$,
which consists of already known dualities for quasi-hereditary algebras).

\begin{corollary}\label{c55.88}
Let $A$ be balanced and $\mathcal{T}(L(i))^{\bullet}$, $i=1,\dots,n$, 
be a complete list of indecomposable tilting objects in $\mathfrak{T}^b$, 
constructed in Theorem~\ref{t55.12}\eqref{t55.12.5}. Then 
$\langle -1\rangle[1]$ induces a (canonical) $\Z$-action on the algebra
\begin{displaymath}
\overline{C}(A)=\End_A\left(\oplus_{l\in\Z}\oplus_{i=1}^n 
\mathcal{T}(L(i))^{\bullet}\langle -l\rangle[l]\right),
\end{displaymath}
which makes $\overline{C}(A)$ into the covering of some algebra $C(A)$.
The algebra $C(A)$ is balanced and  $C(C(A))\cong A$.
\end{corollary}

\begin{proof}
From Theorem~\ref{t55.12} it follows that $C(A)\cong (R(R^{!}))^{opp}$. 
From Lemma~\ref{l55.8} and the assumption that $R(A)^!$ is SCT  it follows that 
the grading on both $R^{!}$ and $C(A)$, induced from $\mathfrak{T}$, is positive. 
In particular, Theorem~\ref{t5.6} and \cite[Theorem~2]{ADL}  now imply that
$C(A)$ is balanced. Since both Ringel and Koszul 
dualities are involutive, we also have $A\cong (R(R(C(A))^{!}))^{opp}$.
\end{proof}

\begin{corollary}\label{c55.89}
Let $A$ be balanced. Then 
$A$ is standard Koszul and $C(A)\cong (A^{!})^{opp}\cong E(A)$.
\end{corollary}

\begin{proof}
$A$ is standard Koszul by Theorem~\ref{t5.6}, in particular, it is
Koszul by \cite[Theorem~1]{ADL}. Further, since no homotopy
is possible in $\mathfrak{T}$, it follows that 
\begin{displaymath}
\ext_A^{l}\left(L(i),L(j)\langle -l\rangle\right)\cong
\Hom_{\mathfrak{T}}\left(\mathcal{T}(L(i))^{\bullet},
\mathcal{T}(L(j))^{\bullet}\langle -l\rangle[l]\right).
\end{displaymath}
The last equality is obviously compatible with the $\Z$-actions and the compositions 
on both sides, which implies that the Koszul dual $(A^{!})^{opp}$ of $A$ is isomorphic 
to $C(A)$.  
\end{proof}

And now we can formulate, probably,  the most surprising 
result of this section.

\begin{corollary}\label{c55.90}
Let $A$ be balanced. Then the algebras $R(A)$,
$E(A)$ and $E(R(A))$ and $R(E(A))$ are also balanced, moreover
\begin{displaymath}
E(R(A))\cong R(E(A))
\end{displaymath}
as quasi-hereditary algebras. In other words, both the Ringel and Koszul 
dualities preserve the class of balanced algebras and commute on this class.
\end{corollary}

\begin{proof}
Follows from Theorem~\ref{t55.12}, Corollary~\ref{c55.88} and 
Corollary~\ref{c55.89}.
\end{proof}

The results, presented in this section motivate the following
natural question: {\em is any SCT=SCK algebra balanced?} 

\section{The graded Ringel dual for the category $\mathcal{O}$}\label{s6} 

In this section we prove that the conditions of Theorem~\ref{t5.5}
are satisfied for the associative algebra, associated with a
block of the BGG category $\mathcal{O}$. To do this we will use the
graded approach to the category $\mathcal{O}$, worked out in
\cite{St}. So, in this section we assume that $A$ is the basic 
associative algebra of an indecomposable integral (not necessarily
regular) block of the BGG category $\mathcal{O}$, \cite{BGG}.
The (not necessarily bijective) indexing set for simple modules will 
be the Weyl group $W$ with the usual Bruhat order (such that the 
identity element is the maximal one and corresponds to the projective 
Verma=standard module). This algebra is Koszul by \cite{BGS,So}, and 
thus we can fix  on $A$ the Koszul grading, which leads us to the 
situation, described in Section~\ref{s5}. Recall that a module, $M$, 
is called {\em rigid} provided that its socle and radical filtrations
coincide, see for example \cite{Ir}. Our main result in this 
section is the following:

\begin{theorem}\label{t6.1}
$\End_A(T)$ is positively graded, moreover, it is generated in
degrees $0$ and $1$. Furthermore, $\nabla$ admits an LT-resolution.
\end{theorem}

\begin{proof}
From \cite[Section~7]{FKM} it follows that 
$T\cong \mathrm{Tr}_{P(w_0)}(P)$ and thus, by Lemma~\ref{l5.1}, 
there is a graded submodule, $T'$ of $P$, which is isomorphic
to $T$ after forgetting the grading. Moreover, again by \cite[Section~7]{FKM},
the restriction from $P$ to $T'$ induces an isomorphism of $\End_A(P)$ and 
$R=\End_A(T)$. So, to prove that $\End_A(T)$ is positively graded it is 
enough to show that $T'\cong T\langle -l\rangle$ for some $l$. Actually,
we will show that this $l$ equals the Loewy length of $\Delta(e)$.

Let $\theta_s$ denote the graded translation functor through the
$s$-wall, see \cite[3.2]{St}. Let $w_0$ denote the longest element
in the Weyl group. The socle of any Verma module in the category
$\mathcal{O}$ is the simple Verma module $\Delta(w_0)$, see 
\cite[Chapter~7]{Di}. This gives, for some $l\in\Z$, a graded 
inclusion, $T(w_0)\langle -l\rangle\cong \Delta(w_0)\langle -l
\rangle\hookrightarrow \Delta(e)$. Moreover, since Verma modules in
$\mathcal{O}$ are rigid by \cite{Ir}, and since their graded filtration
in the Loewy one by \cite[Proposition~2.4.1]{BGS}, it follows that this $l$
equals the Loewy length of $\Delta(e)$. Now we would like to prove by
induction that $T(w_0w)\langle -l\rangle\hookrightarrow P(w)$ for
any $w\in W$. Assume that this is proved for some $w$ and let
$s$ be a simple reflection such that $l(ws)>l(w)$. Translating through
the $s$-wall we obtain
$\theta_s T(w_0w)\langle -l\rangle\hookrightarrow \theta_s P(w)$.
Further, the module  $P(ws)$ is a direct summand of $\theta_s P(w)$ 
(after forgetting the grading). However, from \cite[Theorem~3.6]{St} 
it follows that the inclusion $P(ws)\hookrightarrow\theta_s P(w)$ 
is homogeneous and  has degree $0$. The same
argument implies that the inclusion
$T(w_0ws)\hookrightarrow\theta_s T(w_0w)$ is homogeneous and 
has degree $0$. This gives us the desired inclusion 
$T(w_0ws)\langle -l\rangle\hookrightarrow P(ws)$ of degree $0$
and completes the induction. Adding everything up we obtain a 
graded inclusion of degree $0$ from $T\langle -l\rangle$ to $P$. 

Recall once more that the restriction from $P$ to $T$ induces
an isomorphism of $\End_A(P)$ and $R=\End_A(T)$. Since
$\End_A(P)=A$ is positively graded and is generated in degrees $0$
and $1$, we obtain that $\End_A(T)$ is positively graded and
is generated in degrees $0$ and $1$ as well. 

It is now left to prove the existence of an LT-resolution for 
$\nabla$. Consider the minimal tilting resolution of $\nabla$. 
In Section~\ref{s5} we have defined the grading on $T$ such that 
the canonical projection $T\to \nabla$ 
is a homogeneous map of degree $0$. The kernel of this 
projection is thus graded and has a graded $\nabla$-filtration.
Proceeding by induction we obtain that the minimal tilting resolution
of $\nabla$ is graded. Let $R=R(A)$. Using the functor 
$F=\Hom_A(T,{}_-)$ we transfer this graded tilting resolution 
to a graded projective resolution 
of the direct sum $\Delta^{(R)}$ of standard $R$-modules. By
\cite{So2} we have $A\cong R$, moreover, we have just
proved that the grading on $R$, which is induced from 
$A\mathrm{-gmod}$, is the Koszul one. By \cite[3.11]{BGS}, 
the standard $A$-modules are Koszul, implying that 
the $l$-th term of the projective resolution of 
$\Delta^{(R)}$ is generated in degree $l$. Applying 
$F^{-1}$ we thus obtain an LT-resolution of $\nabla$.
This completes the proof.
\end{proof}

Catharina Stroppel gave an alternative argument for Theorem~\ref{t6.1}
(see Appendix), which uses graded twisting functors. The advantage of 
her approach is that it can be generalized also to the parabolic 
analogue of the category $\mathcal{O}$ defined in \cite{RC}. 

The arguments, used in the proof of Theorem~\ref{t6.1} also imply 
the following technical result:

\begin{corollary}\label{c6.2}
\begin{enumerate}[(1)]
\item The Loewy length $\mathrm{l.l.}(P(w))$ of $P(w)$  
equals $2\mathrm{l.l.}(\Delta(e))-\mathrm{l.l.}(\Delta(w))$. In particular,
for the regular block of $\mathcal{O}$ we have
$\mathrm{l.l.}(P(w))=l(w_0)+l(w)+1$.
\item The Loewy length $\mathrm{l.l.}(T(w))$ of $T(w)$  
equals $2\mathrm{l.l.}(\Delta(w))-1$. In particular,
for the regular block of $\mathcal{O}$ we have
$\mathrm{l.l.}(T(w))=2(l(w_0)-l(w))+1$.
\end{enumerate}
\end{corollary}

\begin{proof}
We start with the second statement. Recall that $\Delta(w)\hookrightarrow 
T(w)$, $T(w)\tto \nabla(w)$, $[T(w):L(w)]=1$ and $L(w)$ is the simple top
of $\Delta(w)$ and the simple socle of $\nabla(w)$. It follows that
$\mathrm{l.l.}(T(w))\geq \mathrm{l.l.}(\Delta(w))+\mathrm{l.l.}(\nabla(w))-1=
2\mathrm{l.l.}(\Delta(w))-1$ since $\mathcal{O}$ has a simple preserving duality.
However, the graded filtration of the tilting module we have just constructed
certainly has semi-simple subquotients (since $A_0$ is positively graded). 
All $\Delta(w')$ occurring in it have 
Loewy length less than or equal to that of $\Delta(w)$ and start in negative 
degrees since $\End_A(T)$ is positively graded. This implies that 
$\mathrm{l.l.}(T(w))\leq 2\mathrm{l.l.}(\Delta(w))-1$ and completes the proof 
of the first part.

Since $P(w)$ has simple top, its graded filtration is the radical one
by \cite[Proposition~2.4.1]{BGS}. However, from the proof of 
Theorem~\ref{t6.1} and from the second part of this corollary, which we have
just proved, it follows that the length of the graded filtration of $P(w)$ 
is exactly  $2\mathrm{l.l.}(\Delta(e))-\mathrm{l.l.}(\Delta(w))$. 

The computations for the regular block follow from the results
of \cite{Ir1} and the proof is complete.
\end{proof}

\begin{corollary}\label{c6.3}
Let $w\in W$. Then the following conditions for $T(w)$ are equivalent:
\begin{enumerate}[(a)]
\item\label{c.6.3.1} $T(w)$ is rigid.
\item\label{c.6.3.2} $\End_A(T(w))$ is commutative.
\item\label{c.6.3.3} $T(w)$ has simple top (or, equivalently, simple socle).
\item\label{c.6.3.35} The center of the universal enveloping algebra surjects 
onto $\End_A(T(w))$.
\item\label{c.6.3.4} $T(w)\hookrightarrow P(w_0)$.
\item\label{c.6.3.5} $P(w_0)\tto T(w)$.
\item\label{c.6.3.6} $[T(w):\Delta(w')]\leq 1$ for all $w'\in W$.
\item\label{c.6.3.7} $[T(w):\nabla(w')]\leq 1$ for all $w'\in W$. 
\item\label{c.6.3.8} $[T(w):\Delta(w_0)]=1$.
\item\label{c.6.3.9} $[T(w):\nabla(w_0)]=1$. 
\end{enumerate}
\end{corollary}

We remark that, though $\Delta(w_0)\cong \nabla(w_0)$ is a simple
module, the numbers $[T(w):\Delta(w_0)]$ and $[T(w):\nabla(w_0)]$ are not the 
composition multiplicities, but the multiplicities in the standard and the 
costandard filtrations of $T(w)$ respectively.

\begin{proof}
By \cite[Section~7]{FKM}, $T(w)\hookrightarrow P(w_0w)$ and the restriction 
induces an isomorphism for the endomorphism rings. Hence the equivalence of
\eqref{c.6.3.2}, \eqref{c.6.3.3}, and \eqref{c.6.3.35} follows from the
self-duality of $T(w)$ and \cite[Theorem~7.1]{St2}. That 
\eqref{c.6.3.3} implies \eqref{c.6.3.1} follows from 
\cite[Proposition~2.4.1]{BGS}. From the proof of Theorem~\ref{t6.1} and
\cite[Theorem~3.6]{St} it follows that the highest and the lowest
graded components of $T(w)$ are one-dimensional. Hence if $T(w)$ does not
have simple top, its graded filtration, which is a Loewy one, does not
coincide with the radical filtration and thus $T(w)$ is not rigid. 
This means that \eqref{c.6.3.1}  implies  \eqref{c.6.3.3}.
Since $L(w_0)$ is the socle of any Verma module, it follows that
\eqref{c.6.3.5} is equivalent to \eqref{c.6.3.3}. And, using the
self-duality of both $T(w)$ and $P(w_0)$ we have that \eqref{c.6.3.5} is 
equivalent to \eqref{c.6.3.4}. 

The equivalence of \eqref{c.6.3.6} and \eqref{c.6.3.7} 
and the equivalence of \eqref{c.6.3.8} and \eqref{c.6.3.9} follows
using the simple preserving duality on $\mathcal{O}$. Since
$[P(w_0):\Delta(w')]=1$ for all $w'$, we get that \eqref{c.6.3.5}
implies \eqref{c.6.3.6}. Let $T(w)$ be such that
\eqref{c.6.3.6} is satisfied. Then, in particular, 
$[T(w):\Delta(w_0)]\leq 1$. Since $L(w_0)$ is a simple socle
of any Verma module, the self-duality of $T(w)$ implies
$[T(w):\Delta(w_0)]=1$, which, in turn, implies that
$T(w)$ has simple top, giving \eqref{c.6.3.3}. Moreover, the same
arguments shows that \eqref{c.6.3.8} implies \eqref{c.6.3.3}.
That \eqref{c.6.3.6} implies \eqref{c.6.3.8} is obvious,
and the proof is complete.
\end{proof}

We remark that (in the case when the equivalent conditions of 
Corollary~\ref{c6.3} are satisfied) the surjection of the center of the 
universal enveloping algebra onto $\End_A(T(w))$ is graded with respect to 
the grading on the center, considered in \cite{So}. 

\begin{corollary}\label{c6.101}
Let $w\in W$, and $s$ be a simple reflection. Then
$\theta_s T(w)=T(w)\langle 1\rangle \oplus
T(w)\langle -1\rangle $ if $l(ws)>l(w)$
and $\theta_s T(w)\in\mathrm{add}(T)$ (as a graded module) otherwise.
\end{corollary}

\begin{proof}
In the case $l(ws)>l(w)$ the statement follows from 
\cite[Section~7]{FKM} and \cite[Section~8]{St2}. If $l(ws)< l(w)$
then Theorem~\ref{t6.1} and \cite[Section~8]{St2} 
implies that $\theta_s T(w)$ has a graded Verma flag, and all
Verma subquotients in this flag are of the form
$\Delta(x)\langle k\rangle$, $k\geq 0$. The self-duality of
$\theta_s T(w)$ now implies that $\theta_s T(w)\in\mathrm{add}(T)$.
\end{proof}

One more corollary of Theorem~\ref{t6.1} is the following:

\begin{proposition}\label{p6.5}
$A$ is a balanced algebra, in particular, all standard, costandard, 
and simple  $A$-modules are $T$-Koszul.
\end{proposition}

\begin{proof}
That standard and costandard $A$-modules are $T$-Koszul follows from 
the fact that $A$ is standard Koszul (see \cite{ADL}) and 
Theorem~\ref{t6.1}. Hence  $A$ is SCT by Theorem~\ref{t6.1} and
Corollary~\ref{c55.89}. Further, the Koszul grading on $A\mathrm{-mod}$ 
induces on $R(A)^!\mathrm{-mod}$ the Koszul grading by 
\cite[Theorem~3]{ADL}. In particular, from Theorem~\ref{t6.1} it follows
that $R(A)^!$ is SCK, that is $A$ is balanced. That simple $A$-modules 
are $T$-Koszul now  follows from Theorem~\ref{t55.12}.
\end{proof}

With the same argument and using the result of Catharina Stroppel presented
in the Appendix, one gets that the algebras of the 
blocks of the parabolic analogue of the category $\mathcal{O}$ in the
sense of \cite{RC} are also balanced. 

We also remark that projective $A$-modules are not $T$-Koszul in general. 
For example, already for $\mathfrak{sl}_2$ we have 
$P(s_{\alpha})\cong T(e)\langle -1\rangle$ and thus $P(s_{\alpha})$ 
is not $T$-Koszul.

\begin{corollary}\label{c6.7}
Let $A$ be the associative algebra of the regular block of the
category $\mathcal{O}$ endowed with Koszul grading. Then the 
category of linear bounded tilting complexes of $A$-modules is
equivalent to $A\mathrm{-gmod}$.
\end{corollary}

\begin{proof}
Since $A$ has a simple preserving duality, it is 
isomorphic to $A^{opp}$, moreover, $A$ is Koszul self-dual
by \cite{So} and Ringel self-dual by \cite{So2}. Hence the
necessary statement follows from Corollary~\ref{c55.7}.
\end{proof}

For singular blocks Corollary~\ref{c55.7} and \cite{BGS}
imply that the category of linear bounded tilting 
complexes of $A$-modules is equivalent to the category of
graded modules over the regular block of the parabolic
category $\mathcal{O}$ with the same stabilizer (and vice versa).

\section{Appendix (written by Catharina Stroppel)}\label{sapp}

In this appendix we reprove Theorem~\ref{t6.1} in a way which implies
the corresponding statement for the parabolic category $\mathcal{O}$. 
Our methods also provide an example for the theory developed in the paper 
in the context of properly stratified algebras. Since we do not use any new
techniques, we refer mainly to the literature. We have to recall several
constructions and definitions.  We restrict ourselves to the case of the
principal block to avoid even more notation. 

For an algebra $A$ we denote by $\mathrm{mod-}A$ ($A\mathrm{-mod-}A$ 
respectively) the category of finitely generated right $A$-modules (finitely 
generated $A$-bimodules). If $A$ is graded, then we denote by $\mathrm{gmod-}A$ 
and $A\mathrm{-gmod-}A$ the corresponding categories of graded modules. 
 
Let $\mathfrak{g}$ be a semisimple Lie algebra with fixed Borel and Cartan 
subalgebras $\mathfrak{b}$, $\mathfrak{h}$, Weyl group $W$ with longest element 
$w_0$, and corresponding category $\mathcal{O}$. Let $\mathcal{O}_0$ be the 
principal block of $\mathcal{O}$ with the simple modules $L(x\cdot0)$ of highest 
weight $x(\rho)-\rho$, where $x\in W$ and $\rho$ denotes the half-sum of positive 
roots. Let $P(x\cdot0)$ be the projective cover of $L(x\cdot0)$. 
  
Let $\mathcal{H}$ denote the category of Harish-Chandra bimodules with generalized
trivial central character from both sides (as considered for example in \cite{SHC}). 
Let $\chi$ denote the trivial central character. For any $n\in\mathbb{Z}_{>0}$ we have
the full subcategories $\mathcal{H}^n$ (and $^n\mathcal{H}$ respectively) of 
$\mathcal{H}$ given by objects $X$ such that $X\mathrm{Ker}\chi^n=0$ ($\mathrm{Ker}\chi^n X=0$ 
respectively). There is an auto-equivalence $\eta$ of $\mathcal{H}$, given by switching 
the left and right action of $U(\mathfrak{g})$ (see \cite[6.3]{Ja}), and giving rise 
to equivalences $\mathcal{H}^n\cong{}^n\mathcal{H}$. For $s$ a simple reflection we have
translation functors through the $s$-wall: $\theta_s$ from the left hand side and 
$\theta_s^r$ from the right hand side (for a definition see \cite[6.33]{Ja} or more 
explicitly \cite[2.1]{Sthom}). In particular, $\eta\theta_s\cong\theta_s^r\eta$. 
Recall the equivalence (\cite[Theorem 5.9]{BG}) $\epsilon:\mathcal{H}^1\cong\mathcal{O}_0$. 
We denote by $L_x=\epsilon^{-1} L(x\cdot 0)$ and consider it also as an object in 
$\mathcal{H}$. Note that $\eta L_x\cong L_{x^{-1}}$ (see \cite[6.34]{Ja}). Let $P^n_x$ 
and ${}^nP_x$ be the projective cover of $L_x$ in $\mathcal{H}^n$ and ${}^n\mathcal{H}$ 
respectively. In particular, $\eta P^n_x\cong {}^n P_x$. 

Recall the structural functor $\mathbb{V}:\mathcal{H}\rightarrow 
S(\mathfrak{h})\mathrm{-mod-}S(\mathfrak{h})$ from \cite{SHC}. We equip the algebra
$S=S(\mathfrak{h})$ with a $\mathbb{Z}$-grading such that $\mathfrak{h}$ is sitting in degree two. 
In \cite{SHC} it is proved that $\mathbb{V} P^n_x$ has a graded lift. By abuse of language, 
we denote the graded lift having $-l(x)$ as its lowest degree also by $\mathbb{V} P^n_x$. Let
$A^n=\mathrm{End}_{S-\mathrm{gmod-}S}(\bigoplus_{x\in W}\mathbb{V} P^n_x)$. Then $A^n$ is a graded
algebra such that $\mathcal{H}^n\cong\mathrm{mod-}A^n$. In particular, $A^1$ is the Koszul
algebra corresponding to $\mathcal{O}_0$ (\cite{BGS}). On the other hand we have
${}^nA=\mathrm{End}_\mathcal{H}(\bigoplus_{x\in W}{}^nP_x)$ and the corresponding equivalence
${}^n\mathcal{H}\cong\mathrm{mod-}{}^nA$.  Concerning the notation we will not
distinguish between objects in $\mathcal{H}^n$ and $\mathrm{mod-}A^n$ or between objects
in $^n\mathcal{H}$ and $\mathrm{mod-}{}^nA$. We fix a grading on $^nA$ such
that $\eta$ lifts to equivalences $\tilde\eta:\mathrm{gmod-}A_n\cong \mathrm{gmod-}A^n$
preserving the degrees in which a simple module is concentrated. More
precisely, $\tilde\eta L(x)\cong L(x^{-1})$, where $L(x)$ denotes the graded
lift of $L_x$, concentrated in degree zero, in the corresponding category. 

Let us fix $n$. For $s$ a simple reflection we denote by $S^s$ the $s$-invariants in $S$. 
We define $\tilde\theta_s:\mathrm{gmod-}A^n\rightarrow\mathrm{gmod-}A^n$ as tensoring 
with the graded $A^n=\mathrm{End}_{S\mathrm{-gmod-}S}(\bigoplus_{x\in W}\mathbb{V} P^n_x)$ bimodule
$\mathrm{Hom}_{S\mathrm{-gmod-}S}(\bigoplus_{x\in W}\mathbb{V} P^n_x,\bigoplus_{x\in W}S\otimes_{S^s}\mathbb{V}
P^n_x\langle -1\rangle)$. Because of \cite[Lemma 10]{SHC}, this is a graded lift 
(in the sense of \cite{St}) of the translation functor 
$\theta_s:\mathcal{H}^n\rightarrow\mathcal{H}^n$. As in \cite{St} we have the adjunction
morphisms $\mathrm{ID}\langle 1\rangle\rightarrow\tilde\theta_s$ and 
$\tilde\theta_s\rightarrow\mathrm{ID}\langle -1\rangle$. Define
$\tilde\theta_s^r=\eta\tilde\theta_s\eta:{}^nA\mathrm{-gmod}\rightarrow
{}^nA\mathrm{-gmod}$. We have again the adjunction morphism $a_s^{(n)}:\mathrm{ID}\langle
1\rangle\rightarrow\tilde\theta_s^r$. Let $T_s^{(n)}$ denote the functor given by
taking the cokernel of $a_s^{(n)}$. We fix a compatible system of surjections
$P^n\tto P^m$ for $n\geq m$. It gives rise to a system of graded projections
$p_{n,m}:{_{}^nA}\tto ^mA$ for $n\geq m$. Let $^\infty A=\varprojlim\;_{}^nA$ and
$^\infty T_s=\varprojlim T^{(n)}:\mathrm{gmod-}^\infty A\rightarrow \mathrm{gmod-}^\infty A$.
Note that $^\infty T_s$ preserves the category $\mathrm{gmod-}A^1$ (considered as a
subcategory of $^\infty A$). In fact, it is a graded lift of Arkhipov's twisting functor 
(as considered in \cite{AS}, \cite{KM}). Let $T_s:\mathrm{gmod-}A^1\rightarrow 
\mathrm{gmod-}A^1$. For $x\in W$ with reduced expression ${[x]}=s_{i_1}s_{i_2}\cdots s_{i_r}$ set
$T_{[x]}=T_{s_1}T_{s_2}\cdots T_{s_r}$. Set $A=A^1$. 

\begin{proposition}\label{prop}
Let $x$, $s\in W$ and let $s$ be a simple reflection. Then the following holds 
\begin{enumerate}
\item The functor $T_{[x]}$ is (up to isomorphism) independent of the chosen
reduced expression.  
\item Moreover, if $sx>x$ and $\Delta(x)\in\mathrm{gmod-}A$ denotes the graded lift
of the Verma module with simple head $L(x)$ (concentrated in degree zero), then
$T_s\Delta(x)\cong\Delta(sx)$ and $T_s\nabla(sx)\cong\nabla(x)$, where
$\nabla(x)$ denotes the graded lift of the dual Verma module with socle
$L(x)$ (concentrated in degree zero). 
\end{enumerate}
\end{proposition}

\begin{proof}
We consider now the adjunction morphism $b_s:\mathrm{ID}\rightarrow\tilde\theta_s^r$ 
between endofunctors on $\mathrm{mod-}_{}^nA$. Let $\tilde{T}_s$ denote the functor given by
taking the cokernel of $b_s$, restricted to $\mathrm{mod-}A$. Let $\tilde T_{[x]}=\tilde{T}_{s_1}\tilde{T}_{s_2}\cdots\tilde{T}_{s_r}$. Then $\tilde{T}_{[x]}$ does 
not depend on the chosen reduced expression (\cite{Joseph}, \cite{KM}). If we show that 
$\tilde{T}_{[x]}$ is indecomposable, then a graded lift is unique up to isomorphism and 
grading shift, and the statement follows say from the second part of the theorem. Set 
$G=\tilde{T}_x$. Let us prove the indecomposability: We claim that the canonical evaluation 
morphism $\mathrm{End}(G)\rightarrow \mathrm{End}_\mathfrak{g}(G P(w_0\cdot0))$, 
$\phi\mapsto\phi_{P(w_0\cdot0)}$, is an isomorphism. Assume $\phi_{P(w_0\cdot0)}=0$. Let $P$ 
be a projective object in $\mathcal{O}_0$. Then there is a short exact sequence 
\begin{equation}\label{eq:ses}
P\rightarrow\oplus_{I} P(w_0\cdot0)\rightarrow Q
\end{equation}
for some finite set $I$ and some module $Q$ having a Verma flag. (To see this consider the 
projective Verma module. It is the unique Verma submodule of ${P(w_0\cdot0)}=0$, hence the desired 
sequence exists. The existence of the sequence for any projective object follows then using 
translation functors.) By \cite[Lemma 2.1]{AS}, we get an exact sequence  $G P\rightarrow\oplus_{I}
GP(w_0\cdot0)\rightarrow GQ$. Hence  $\phi_{P(w_0\cdot0)}=0$ implies $\phi_P=0$ for any 
projective object $P$. Since $G$ is right exact, it follows  $\phi=0$. Let 
$g\in\mathrm{End}_\mathfrak{g}(G P(w_0\cdot0))$. Since $\mathrm{End}_\mathfrak{g}(G P(w_0\cdot0))\cong\mathrm{End}_\mathfrak{g}(P(w_0\cdot0))$ (\cite[Proposition 5.3]{AS}), $g$ defines
an endomorphism of $G$ when restricted to the additive category generated by $P(w_0\cdot0)$. Note 
that (by taking the injective hull of $Q$) the sequence  \eqref{eq:ses} gives rise to an exact sequence 
\begin{equation*}
0\rightarrow P\rightarrow\oplus_{I} P(w_0\cdot0)\rightarrow \oplus_{I'} P(w_0\cdot0)
\end{equation*}    
for some finite sets $I$, $I'$. Using again \cite[Lemma 2.1]{AS} we get an exact sequence 
\begin{equation*}
0\rightarrow G P\rightarrow\oplus_{I} G P(w_0\cdot0)\rightarrow\oplus_{I'} G P(w_0\cdot0).
\end{equation*}  
Hence $g$ defines an endomorphism $g_P$ of $P$. Standard arguments show that this is 
independent of the chosen exact sequence. Since $G$ is right exact, $g$ extends uniquely 
to an endomorphism $\phi$ of $G$. By construction $\phi_{P(w_0\cdot0)}=g$. This proves 
the surjectivity. Since $\mathrm{End}_\mathfrak{g}(G P(w_0\cdot0))\cong
\mathrm{End}_\mathfrak{g}(P(w_0\cdot0))$ is a local ring, the functor $G$ is indecomposable. 
This proves the first part of the proposition. 

We have $\tilde{T}_s f(\Delta(x))\cong f(\Delta(sx))$, where $f$ denotes the grading forgetting 
functor. Hence, $T_s(\Delta(x))\cong \Delta(sx)\langle k\rangle$ for some $k\in\mathbb{Z}$. On the
other hand $\eta T_s\eta\Delta(x^{-1})\cong\Delta(x^{-1}s)$ (\cite[Theorem~3.6]{St}). Hence 
$k=0$ and $T_s\Delta(x)\cong\Delta(sx)$. Forgetting the grading we have $\tilde{T}_s f(\nabla(sx))\cong
f(\nabla(x))$. On the other hand $\eta T_s\eta\nabla((sx)^{-1})\cong\nabla(x^{-1})$ 
(\cite[Theorem 3.10]{St}). The second part of the proposition follows.  
\end{proof}

Since $T_{[x]}$ does not depend on the chosen reduced expression, we denote it
just $T_x$ in the following. Let $P(x)\in\mathrm{gmod-}A$ be the indecomposable
projective module with simple head $L(x)$ concentrated in degree zero. Set
$P=\bigoplus_{x\in W}P(x)$. Let $T(x)$ denote the graded lift of an
indecomposable tilting module, characterized by the property that $\Delta(x)$
is a submodule and $\nabla(x)$ is a quotient. Let $T=\bigoplus_{x\in W}
T(x)$. 

\begin{theorem}\label{tapp}
Let $x\in W$. There is an isomorphism of graded algebras
\begin{eqnarray*}
\mathrm{End}_{A}(P)\cong\mathrm{End}_{A}(T_x P).
\end{eqnarray*}
For $x=w_0$ we get in particular 
\begin{eqnarray*}
\mathrm{End}_{A}(P)\cong\mathrm{End}_{A}(T).
\end{eqnarray*}
\end{theorem}

\begin{proof}
The first isomorphism follows directly from \cite[Lemma 2.1]{AS} and the definition of
$T_x$. For the second we claim that $T_{w_0} P(y)\cong T(w_0y)$. By 
Proposition~\ref{prop} we have $T_{w_0}P(0)\cong \Delta(w_0)$. Hence, the statement 
is true for $y=e$. Using translation functors we directly get $T_{w_0} P(y)\cong
T(w_0y)\langle k\rangle$ for some $k\in\mathbb{Z}$. On the other hand $P(y)$ surjects onto 
$\Delta(y)$. Then $T_{w_0}P(y)$ surjects onto $T_{w_0}\Delta(y)$. The latter is 
isomorphic to  $\nabla(w_0y)$, since  $\Delta(w_0)\cong\nabla(w_0)$. 
\end{proof}

Let $\ppp$ be a parabolic subalgebra of $\mathfrak{g}$ with corresponding parabolic
subgroup $W_\ppp$ of $W$. Let $\mathcal{O}_0^\ppp$ be the full subcategory of $\mathcal{O}_0$
given by locally $\ppp$-finite objects. If $P\in\mathcal{O}_0$ is a minimal projective
generator, then its maximal quotient $P^\ppp$ contained in $\mathcal{O}_0^\ppp$ is a
minimal projective generator of $\mathcal{O}_0^\ppp$ and $\mathrm{End}_\mathfrak{g}(P^\ppp)$ 
inherits a grading from $A=\End_\mathfrak{g}(P)$. We will consider then the category
$\mathrm{gmod-}A^\ppp$ as the full subcategory of $\mathrm{gmod-}A$ given by all objects having
only composition of the form $L(x)\langle k\rangle$, where $k\in\mathbb{Z}$ and $x\in W^\ppp$, the 
set a shortest coset representative of $W_\ppp\backslash W$. Let $\Delta^\ppp(x)\in\mathrm{gmod-}A^\ppp$,
$\nabla^\ppp(x)$ be the standard graded lifts of the standard and costandard modules in 
$\mathcal{O}_0^\ppp$ (which were denoted by $\Delta(x)$ and $\nabla(x)$ in Section~\ref{s6}). Let 
$T^\ppp$ be the module $T$ from Corollary~\ref{c6.101} for the category $\mathrm{gmod}-A^\ppp$.   
Then Theorem~\ref{t6.1} generalizes to the following   

\begin{corollary}\label{capp}
$\mathrm{End}_{A^\ppp}(T^\ppp)$ is positively graded, moreover, it is generated in degrees
$0$ and $1$. Furthermore, $\nabla$ admits an LT-resolution.  
\end{corollary}

\begin{proof}
Let $w=w_0^\ppp\in W^\ppp$ be the longest element. Then $\Delta^\ppp(w)$
is a tilting module and canonically a quotient of $\Delta(w)\cong
T_w\Delta(e)=T_wP(e)$. Using translation functors we get that $T^\ppp$ is a
quotient of $T_w P$. Hence, there is a surjection of graded algebras from
$\mathrm{End}_{A}(T_wP)\cong\mathrm{End}_{A}(P)$ onto
$\mathrm{End}_{A}(T^\ppp)$. Hence $\mathrm{End}(T_w P^\ppp)\cong\End_{A^\ppp}(T^\ppp)$
is positively graded and generated in degrees 0 and 1. The existence of the
resolution follows using the same arguments as in the proof of Theorem~\ref{t6.1}.
\end{proof}

\vspace{1cm}

\begin{center}
{\bf Acknowledgments}
\end{center}

The research was done during the visit of the second author to Uppsala
University, which was partially supported by the Royal Swedish 
Academy of Sciences, and by The Swedish Foundation for International
Cooperation in Research and Higher Education (STINT).  This support and 
the hospitality of Uppsala University are gratefully acknowledged. The 
first author was also partially supported by the Swedish Research Council. 
We thank Catharina Stroppel for many useful remarks and comments on the
preliminary version of the paper and for writing the Appendix.

\vspace{0.5cm}

\noindent
Volodymyr Mazorchuk, Department of Mathematics, Uppsala University,
Box 480, 751 06, Uppsala, SWEDEN, 
e-mail: {\tt mazor\symbol{64}math.uu.se},
web: {``http://www.math.uu.se/$\tilde{\hspace{1mm}}$mazor/''}.
\vspace{0.2cm}

\noindent
Serge Ovsienko, Department of Mechanics and Mathematics, Kyiv Taras
Shevchenko University, 64, Volodymyrska st., 01033, Kyiv, Ukraine,
e-mail: {\tt ovsienko\symbol{64}zeos.net}.

\end{document}